\documentclass[a4paper, 12pt, reqno]{amsart}

\usepackage[protrusion=false,final]{microtype}
\usepackage[margin=1in]{geometry}
\usepackage{graphicx}
\usepackage{xcolor}
\usepackage{parskip}
\usepackage{booktabs} 
\usepackage{amssymb}
\usepackage{amsmath}
\usepackage{mathtools}
\usepackage{hyperref}
\usepackage{cleveref}
\usepackage{enumerate}
\usepackage{enumitem}
\usepackage{bbm}

\numberwithin{equation}{section}

\theoremstyle{plain}
\begingroup
\newtheorem{theorem}{Theorem}[section]
\newtheorem{lemma}[theorem]{Lemma}
\newtheorem{proposition}[theorem]{Proposition}
\newtheorem{corollary}[theorem]{Corollary}
\endgroup
\theoremstyle{definition}
\begingroup
\newtheorem{definition}[theorem]{Definition}

\endgroup
\theoremstyle{remark}
\begingroup
\newtheorem{remark}[theorem]{Remark}
\endgroup

\renewcommand{\o}{\Omega}

\newcommand\e{\varepsilon}
\newcommand{\R}{\mathbb{R}}
\newcommand{\N}{\mathbb{N}}
\newcommand{\Z}{\mathbb{Z}}

\newcommand\iO{\int_\Omega}
\newcommand{\iOQ}{\int_{\Omega}\int_Q}

\newcommand{\hno}{{\mathcal H}^{N-1}}

\newcommand{\wkstr}{\xrightharpoonup{*}}
\newcommand{\wk}{\rightharpoonup}

\newcommand{\liminfn}{\underset{n \to \infty}{\liminf}\;}
\newcommand{\limsupn}{\underset{n \to \infty}{\limsup}\;}
\newcommand{\limn}{\underset{n \to \infty}{\lim}\;}

\newcommand{\unf}{\mathcal{U}_{\delta_n}}
\newcommand{\tr}{\mathcal{T}}
\newcommand{\oh}{{\hat{\o}_\delta}}
\newcommand{\ohn}{{\hat{\o}_{\delta_n}}}
\DeclarePairedDelimiter\floor{\lfloor}{\rfloor}
\DeclareMathOperator*{\esssup}{ess\,sup}

\title[The supercritical case - Fixed wells]{Homogenization and phase separation with fixed wells - The supercritical case}
\author[R. Cristoferi] {Riccardo Cristoferi} 
\address[R. Cristoferi]{\\Department of Mathematics - IMAPP\\ Radboud University\\
Heyendaalseweg \\ 6525 AJ Nijmegen\\ The Netherlands}
\author[I. Fonseca] {Irene Fonseca} 
\address[I. Fonseca]{Department of Mathematics\\ Carnegie Mellon University\\
Forbes Avenue\\Pittsburgh PA 15213, USA}
\author[L. Ganedi] {Likhit Ganedi} 
\address[L. Ganedi]{Institut für Mathematik\\ RWTH Aachen\\
Templergraben\\52062 Aachen, Germany}

\keywords{homogenization, phase separation, $\Gamma$-expansion}

\begin{document}

	\maketitle
 
\begin{abstract}
    A variational model for the interaction between homogenization and phase separation is considered in the regime where the former happens at a finer scale than the latter. The first order $\Gamma-$limit is proven to exhibit a separation of scales which has been previously conjectured in \cite{BradiesZeppieri}\cite{Hagerty}.
\end{abstract}


\section{Introduction}
Composite materials are important to modern technology as the mixing of two different material properties at fine scales can give rise to unexpected emergent behavior \cite{clyne2019introduction}. Therefore, understanding the process of phase separation on such materials is crucial to leveraging these processes for technological applications. 

For a homogeneous material, the distribution of stable phases is commonly modeled by using the Cahn-Hilliard free energy (also known as the Modica-Mortola functional, in the mathematical community). The energy reads as
\[
E_{\e}(u) \coloneqq \int_\o{\left[W\left(u(x)\right)+\e^2|\nabla u(x)|^2\right]\;dx},
\]
where $u\in W^{1,2}(\Omega;\R^M)$ represents the distribution of phases, $\varepsilon>0$ is a small parameter that represents the width of the transition layers between the different materials that form the composite, and the free energy $W:\R^M\to[0,\infty)$ that vanishes at the stable critical phases.
It was first proved by Modica and Mortola \cite{modica87}\cite{modica77} in the scalar case that in the limit this energy minimizes perimeter, i.e, interfacial energy. This sharp interface limit was conjectured by Gurtin \cite{Gurtin} to hold in more generality, and was later proven in \cite{KohnStern, Stern, FonsecaTartar}. Since then, many variants have been studied such as having multiple phases \cite{baldo}, fully coupled singular perturbations \cite{Bou, FonPopo}, and even the case in which the wells of $W$ are allowed to depend on position \cite{CriGra,CristoferiFonsecaGanedi}.

Our interest here is in a heterogeneous material where the heterogeneities are modeled with an oscillating periodic potential $W:\R^N\times\R^M\to[0,+\infty)$ that is $Q$-periodic in the first variable, where $Q\subset\R^N$ is the unit cube $(-\frac{1}{2},\frac{1}{2})^N$. We consider the wells of $W$ to be fixed constants $a,b \in \R^M$, while the case in which the  wells of $W$ are dependent on the spatial variable will be considered in future work.

The energy functional reads as
 \[
 \mathcal{F}_{\e,\delta}(u) \coloneqq \int_\o\left[W\left(\frac{x}{\delta},u\right)+\e^2|\nabla u|^2\right]\;dx.
 \]

In order to analyze the behavior of minimizers, we consider a $\Gamma-$expansion \cite{braides2008asymptotic,anzellotti1993asymptotic}. At the lowest order, it is possible to show that, $F_{\e,\delta} \xrightarrow{\Gamma} F_0$ as $\delta,\e\to0$ with 
$$F_0(u) := \int_\o{W_r(u)\;dx},$$
where $u \in L^1(\o;\R^M)$ and $W_r$ is a homogenized potential whose form depends on the rate of convergence $r: = \lim \frac{\delta}{\e}$. Since we are in the regime of fixed wells $a,b$, it is always possible to find many minimizers of $F_0$ which achieve zero energy even with a mass constraint. Thus, in order to better understand the minimizers, we need to consider the next order in the $\Gamma-$expansion. Similar to the heuristics for the homogeneous Modica-Mortola functional, ìt is possible to determine that the energy of having a transition layer between the phases will be of order $\e$.
This leads us to consider the rescaled functional
\begin{equation}{\label{rescaled_fnal}}
    {\mathcal{E}}_{\varepsilon,\delta}(u) \coloneqq \frac{1}{\varepsilon}\mathcal{F}_{\e,\delta}(u)
 = \int_\o \left[\frac{1}{\varepsilon}W\left(\frac{x}{\delta},u\right)
 +\e|\nabla u|^2 \right]\;dx
\end{equation}
However, this energy has not been studied much in the literature due to technical mathematical difficulties it poses. The behavior of minimizers of such a model depends greatly on the rate $r$ at which $\e,\delta$ comparatively decay to $0$, i.e., $\e\ll \delta,\e \sim \delta$, or $\delta \ll \e$. In \cite{CriFonHagPop} the first two authors in collaboration with Hagerty and Popovici, rigorously characterized the first order $\Gamma$-limit when $\e \sim \delta$. This has been recently extended to the fully coupled scalar case with stochastic homogenization \cite{Mar22}. For $\e \ll \delta$, the characterization of the $\Gamma-$limit is still open, but in \cite{CristoferiFonsecaGanedi}, the authors have identified an intermediate scaling of $\frac\e\delta$ and characterized the $\Gamma-$limit with respect to strong two-scale convergence in an analysis that also extends to the case of spatially dependent phases.

In this paper, we study the case $\delta\ll\varepsilon$, and we prove a separation of scales that has only been conjectured thus far (see \cite{BradiesZeppieri}\cite{Hagerty}), namely that the first order $\Gamma$-limit is the $\Gamma$-limit of the functional
\[
u \mapsto \int_\o{ \left[\frac{1}{\varepsilon} W_{\mathrm{hom}}(u(x))
    + \e|\nabla u(x)|^2 \right]\;dx},
\]
where $W_{\mathrm{hom}}$ is the homogenized potential of $W$ defined as
\[
W_{\mathrm{hom}}(z)\coloneqq \int_Q W(y,z) dy.
\]
Heuristically, this is expected because the regime $\delta\ll\varepsilon$ suggests that we first homogenize (namely, we first send $\delta\to0$), and then we study phase separation (namely, we send $\varepsilon\to0$). Indeed, \cite{Hagerty} was able to use the technique of direct replacement of the potential by $W_\mathrm{hom}$ (first used in \cite{BradiesZeppieri} in a similar setting) to show the $\Gamma-$limit when $\delta \ll \e^\frac{3}{2}$. Here, we are able to prove that the same heuristics can be made rigorous even when $\delta\ll\varepsilon$, by using an intermediate cell problem derived from the two-scale unfolding of the functional.

Our strategy, which will be outlined in \Cref{sec:Strategy}, enables us to weaken the regularity requirements of $W$ to be only Carath\'eodory and to remove restrictions such as quadratic behavior near the wells, which are sometimes used in literature \cite{CristoferiFonsecaGanedi}\cite{FonsecaTartar}. Another important feature is that we prove strong compactness directly from the coercivity and polynomial growth bounds rather than assuming the existence of a continuous double well potential $W_H$ independent of $x$ and such that $W_{H}(z)\leq W(x,z)$ for every $z \in \R^M$, as used in \cite{Hagerty}.

Finally, we note that sometimes in the literature the heterogeneity is entered into the energy through the singular perturbation \cite{AnsBraChi1}\cite{AnsBraChi2}. This creates a similar separation of scales effect, but it requires different techniques and leads to anisotropic effective limits. The techniques used in \cite{AnsBraChi1}{\cite{AnsBraChi2}  has since also been applied to the case of the Ambrosio-Totorelli energy \cite{Bach2022}\cite{BacMarZep22}.


\subsection{Main results}

Let $\o\subset\R^N$ be a bounded open set, and $N, M \geq 1$.
Denote by $Q\coloneqq (-1/2,1/2)^N$ the unit cube in $\R^N$ centered at the origin.
Let $W:\R^N\times\R^M\to[0,+\infty)$ be a measurable function which satisfying the following hypotheses:
	\begin{enumerate}[label=(W\arabic*), ref=(W\arabic*)]
            \item \label{W1} $W$ is a Carath\'eodory function which is $Q$-periodic in the spatial variable, \emph{i.e.},
            \begin{itemize}
                \item $z\mapsto W(x,z)$ is continuous for $\mathcal{L}^N-$a.e. $x\in Q$,
                \item $x\mapsto W(x,z)$ is measurable and $Q-$periodic for all $z \in \R^M$.
            \end{itemize}
            \item \label{W2} There are $a,b \in \R^M$ such that 
            \[
            W(x,z)=0 \iff z \in \{a,b\}.
            \]
		\item \label{W3} There exists $R>0$ such that for $\mathcal{L}^N$-a.e. $x\in Q$, 
  \[
W(x,z)\geq \frac{1}{R}|z|,
\]
if $|z|\geq R $.
\item\label{W4} For every $M>0$, there exists a constant $C_M > 0$ depending only on $M$ such that 
$$\esssup_{x \in Q, |z|\leq M} W(x,z)\leq C_M$$ 
\end{enumerate}

\begin{remark}
We note that we use two-wells for convenience. For multiple wells, a similar result holds using \cite{baldo}.
\end{remark}

\begin{remark}
Assumption \ref{W1} is the minimum regularity requirement when allowing for material inclusions which are a common form of composite materials.
The coercivity growth assumption \ref{W3} is standard in such type of problems. We note that \ref{W4} allows for a general class of upper bounds on $W$ including polynomial growth conditions.
Note that we are not assuming any behavior of the potential close to the wells. This is a novelty in comparison with previous works (see,e.g., \cite{FonsecaTartar}\cite{CriFonHagPop}\cite{CristoferiFonsecaGanedi}).
\end{remark}

We now introduce the functionals that we will study.

\begin{definition}\label{def:functional_En}
Let $\{\e_n\}_{n}, \{\delta_n\}_{n}$ be infinitesimal sequences such that
\[
\lim_{n\to\infty}\frac{\delta_n}{\e_n}= 0.
\]
For $n\in\N$, define the functional $\mathcal{E}_n: L^1(\o;\R^M)\to[0,+\infty]$ as
\[
\mathcal{E}_n (u) \coloneqq \begin{dcases} 
\int_\o \left[\, \frac{1}{\e_n} W\left( \frac{x}{\delta_n}, u(x) \right)  \,+ \e_n|\nabla u(x)|^2 \,\right] dx & \textup{if } u\in W^{1,2}(\o;\R^M),\\[5pt]
+\infty & \textup{else.}
\end{dcases}
\]
\end{definition}

We next introduce the limiting functional in \eqref{eq:limiting_functional}.

\begin{definition}\label{def:limiting_constant}
For $z\in\R^M$, let
\[
{W}_{\mathrm{hom}}(z)\coloneqq \int_Q W(x,z)\, dx.
\]
Set
\[
\sigma_{\mathrm{hom}}\coloneqq \inf\left\{\int_{-1}^1 2\sqrt{{W}_{\mathrm{hom}}(\gamma(t))}|\gamma'(t)|dt
    : \gamma\in \mathrm{Lip}_{\mathcal{Z}}([-1,1];\R^M),\, \gamma(-1)=a,\, \gamma(1)=b \right\},
\]
where $\mathrm{Lip}_{\mathcal{Z}}([-1,1];\R^M)$ is the space of continuous curves $\gamma:[-1,1]\to\R^M$ such that $\gamma\in \mathrm{Lip}(T,\R^M)$, for every compact set $T\subset[-1,1]$ disjoint from $\{t\in[-1,1] : \gamma(t)\in\{a,b\}\}$.
\end{definition}

\begin{definition}\label{def:limiting_energy}
Define the functional $\mathcal{E}_\infty: L^1(\o;\R^M)\to[0,+\infty]$ as
\begin{equation}\label{eq:limiting_functional}
\mathcal{E}_\infty (u) \coloneqq \begin{dcases} 
\sigma_{\mathrm{hom}} \text{Per}(\{u = a\};\o) & \textup{if } u\in BV(\o;\{a,b\}),\\[5pt]
+\infty & \textup{else,}
\end{dcases}
\end{equation}
where $\text{Per}(\{u = a\};\o)$ denotes the perimeter of the set $\{u = a\}$ in $\o$.
\end{definition}

We are now in position to state the two  main results of this paper, namely pre-compactness of sequences with uniformly bounded energy, and the $\Gamma$-convergence of $\{\mathcal{E}_n\}_n$.

\begin{theorem}\label{thm:main_compactness}
Let $\{\e_n\}_{n}, \{\delta_n\}_{n}$ be infinitesimal sequences such that
\[
\lim_{n\to\infty}\frac{\delta_n}{\e_n}= 0.
\]
Let $\{u_n\}_n\subset L^1(\o;\R^M)$ be such that
\[
\sup_{n\in\N} \mathcal{E}_n(u_n) < \infty.
\]
Then, there exist $u\in BV(\Omega;\{a,b\})$ and a subsequence $\{u_{n_k}\}_k$ such that $u_{n_k}\to u$ strongly in $L^1(\o;\R^M)$.
\end{theorem}

\begin{theorem}\label{thm:main_gamma}
Let $\{\e_n\}_{n}, \{\delta_n\}_{n}$ be infinitesimal sequences such that
\[
\lim_{n\to\infty}\frac{\delta_n}{\e_n}= 0.
\]
Then, $\mathcal{E}_n\stackrel{\Gamma}{\rightarrow} \mathcal{E}_\infty$ with respect to the strong  $L^1(\o;\R^M)$ convergence.
\end{theorem}

The strategy of the proofs are stable enough to allow for a mass constraint to be incorporated in the functional.

\begin{definition}\label{def:mass_constrain_functional}
 Let $m\in(0,|\Omega|)$. We define the mass constrained functional
\[
\widehat{\mathcal{E}}_n (u;m) \coloneqq \begin{dcases} 
\mathcal{E}_n(u) & \textup{if } u\in W^{1,2}(\o;\R^M),\, \int_\Omega u\, dx = ma+(1-m)b,\\[5pt]
+\infty & \textup{else.}
\end{dcases}
\]

We also define the suitable limiting problem
\[
\widehat{\mathcal{E}}_\infty (u;m) \coloneqq \begin{dcases} 
\mathcal{E}_\infty(u),  & \textup{if } u\in BV(\o;\{a,b\}),\, \int_\Omega u\, dx = ma+(1-m)b,\\[5pt]
+\infty & \textup{else.}
\end{dcases}
\]   
\end{definition}
\begin{corollary}\label{cor:mass_cosntrain}
Let $m\in(0,|\Omega|)$ 
and $\{\e_n\}_{n}, \{\delta_n\}_{n}$ be infinitesimal sequences such that
\[
\lim_{n\to\infty}\frac{\delta_n}{\e_n}= 0.
\]
Then, $\widehat{\mathcal{E}}_n\stackrel{\Gamma}{\rightarrow} \widehat{\mathcal{E}}_\infty$ with respect to strong  $L^1(\o;\R^M)$ convergence
Moreover, pre-compactness for sequences with uniformly bounded $\widehat{\mathcal{E}}_n$ energy holds.
\end{corollary}
\begin{remark}
The $\Gamma$-convergence results stated in Theorem \ref{thm:main_gamma} and Corollary \ref{cor:mass_cosntrain} allow to get the standard convergence of minima and minimizers (see \cite[Corollary 7.20]{Dalmasobook})), as well as approximation of isolated local minimizers (see \cite{KohnStern}).
\end{remark}


 \subsection{Outline of the Strategy}\label{sec:Strategy}
    The recovery sequence is the same recovery sequence as for the Modica-Mortola energy with potential $W_\mathrm{hom}$. This, and the modifications required to satisfy the usual mass constraint are detailed in Sections 6 and 7, which require some care due to the minimal assumptions on $W$. The core of the work is in proving the Liminf inequality in Section 4 (see \Cref{thm:liminf}), with some useful preliminary and auxiliary results contained in Sections 2 and 3. Here, we outline the main ideas of the proof of the Liminf inequality.
    
	Take a sequence $\{u_n\} \subset W^{1,2}(\Omega;\R^M)\cap L^\infty(\o;\R^M)$ with bounded energy which achieves the $\Gamma$-Liminf (see \Cref{def:gc}). We can partially unfold the energy with the unfolding operator (see \Cref{def:unfold}) just on the potential, and use the non-negativity to throw away the boundary terms at the cost of the correct inequality,
 
	\begin{equation*}\label{eq:partial_unf}
		\mathcal{E}_n[u_n]\geq \int_\o{ \left[\int_Q{\frac{W\left(y,\unf u_n\right)}{\e_n}\;dy}+\e_n|\nabla u_n|^2 \right]\;dx}.
	\end{equation*}
	
	Applying Young's inequality, we get
	$$\mathcal{E}_n[u_n]\geq \int_\o{2\left[\int_Q W(y,\unf u_n(x,y))\;dy\right]^{\frac{1}{2}}|\nabla u_n|\;dx}=: F_n[u_n].$$
	
	To finish we would need to replace the integral under the square root by $$W_{\mathrm{hom}}(u_n) = \int_Q W(y,u_n)\;dy.$$ 

    First, we notice that the following term is negligible in the limit (see \eqref{eq:key_new_estimate}),  
	$$\int_\o{\int_Q |\unf u_n-u_n|\;dy |\nabla u_n|\;dx} \to 0,$$
    and we rewrite $F_n$ as
		$$\int_\o{2\left[\int_Q W(y,u_n(x)+\unf u_n(x,y)-u_n(x))\;dy\right]^{\frac{1}{2}}|\nabla u_n|\;dx}.$$
	The idea is to claim that the potential essentially acts like $W_\mathrm{hom}(u_n)$, but with the exception of some small sets. We fix an $\eta>0$. Using a slicing argument, we find a sequence $\{v^\eta_n\}_n \subset L^\infty(\Omega;W^{1,2}_0(Q))$ with $\|v_n\|_\infty \leq \eta $ and with lower energy,
	 $$\liminfn F_n[u_n]\geq \liminfn\int_\o{2\left[\int_Q W(y,u_n+v_n^\eta)\;dy\right]^{\frac{1}{2}}|\nabla u_n|\;dx}.$$
	We define a double-well function $W^\eta(z)$ such that $v_n^\eta$ is admissible in the infimum (see Section 3.2). After appropriate truncations and removal of small "bad sets", we achieve 
    $$\liminfn F_n[u_n]\geq \liminfn\int_\o{2\sqrt{ W^\eta(u_n)}|\nabla u_n|\;dx}.$$
    Since $W^\eta$ is still a continuous double-well function with the same wells (see \Cref{thm:prop_W_eta}), we can apply "classical" compactness and $\Gamma-$ Liminf arguments \cite{FonsecaTartar} to obtain that up to a subsequence, $u_n \to u$ strongly in $L^1(\o;\R^M)$ and   
	 $$\liminfn F_n[u_n] \geq \sigma_\eta \text{Per}(\{u = a\};\o).$$
  
	To conclude, we show that $\sigma_\eta \nearrow \sigma_{\mathrm{hom}}$, where $\sigma_{\mathrm{hom}}$ is as defined in \Cref{def:limiting_energy}.
	
    Here we use ideas from \cite{SternbergZuniga}, where Zuniga and Sternberg showed under very minimal conditions the existence of a minimizer of the geodesic problems that underly $\sigma_\eta$. We recall these properties in Section 2.3, and use them to prove some critical results in Section 3.3.  
  
	In the sequel, we will often take subsequences without relabeling, $C$ will be a generic constant that may change between inequalities, and subscripts to $C$ will describe the limiting parameters that it depends upon.


\section{Preliminaries}

\subsection{The unfolding operator}

We recall the unfolding operator, which was first used to pass to the limit in periodic homogenization problems by rendering the microscopic scale to behave macroscopically. Two scale convergence was shown to be characterized as $L^p$ convergence on the product space through the unfolding operator \cite{CioDamGri08}\cite{Vis06}\cite{Vis07}. While we do not need two scale convergence in this paper, the unfolding operator provides here a useful tool to encode the usual change of variable used in homogenization problems. 

\begin{definition}\label{def:unfold}
For $\delta >0$, let
\[
\oh\coloneqq\bigcup_{z_i\in I_\delta} \left(\overline{z_i+\delta Q}\right)
\]
where $I_\delta$ is the set of points $k\in\delta\Z^N$ such that $\overline{k+\delta Q}\subset\o$.
The unfolding operator $\mathcal{U}_{\delta}:L^1(\o;\R^M)\to L^1(\o;L^1(Q;\R^M))$ is defined as
\begin{equation}\label{eq:unfolding-operator}
    \mathcal{U}_{\delta}(u)(x,y)\coloneqq\left\{
\begin{array}{ll}
u\Big(\delta\floor[\Big]{\frac{x}{\delta}}+\delta y\Big) & \text
{ for }x\in\oh,\; y\in Q, \\
& \\
a & \text{ if } x\in\o\setminus\oh, y\in Q,
\end{array}
\right.
\end{equation}
where, given an enumeration $\{k_i\}_{i\in\N}$ of $\Z^N$,
\begin{equation}\label{eq:floor}
\floor{x}\coloneqq k_i\,\quad\quad\quad i\coloneqq \min\left\{j\in\N \,:\, k_j\in \mathrm{argmin}\{ |k-x| : k\in\Z^N \}\,\right\}
\end{equation}
is the integer part of $x\in\R^N$, and $a$ is the well in \ref{W2}.
\end{definition}

\begin{remark}\label{rmk:unfolding}
This definition of the unfolding operator is nonstandard as we make the unfolding operator nonzero in the small boundary set $\Lambda_\delta \times Q$. This has been used previously in \cite{CristoferiFonsecaGanedi} to simplify some of the computations. In particular, by change of variables and periodicity of W we can rewrite
\begin{align*}
   \iO{W\left(\frac{x}{\delta},u\right)}\;dx&\geq \int_\oh{W\left(\frac{x}{\delta},u\right)}\;dx = \sum_{z_i \in I_\delta} \delta^N \int_{Q}{W\left(y,u(z_i + \delta y)\right)\;dy}\\
   &= \sum_{z_i \in I_\delta}\int_{z_i + \delta Q}\int_Q{W\left(y,u(z_i + \delta y)\right)\;dydx} = \int_\oh\int_Q{W\left(y,\mathcal{U}_{\delta} u\right)\;dydx} \\
   &= \iOQ{W\left(y,\mathcal{U}_{\delta} u\right)\;dydx},  
\end{align*}
where in the second to last equality we used that for every $x \in z_i + \delta Q$, $\delta\floor[\Big]{\frac{x}{\delta}}= z_i$.
In the last equality, we are able to add back in the boundary set to the unfolded integral as by definition of our unfolding operator, we have $W(y,\mathcal{U}_{\delta} u)= W(y,a) = 0$.

Furthermore, we note that for Sobolev functions, gradients transform by chain rule to be
$$\delta \mathcal{U}_\delta(\nabla u(x)) = \nabla_y(\mathcal{U}_\delta u(x,y)).$$
\end{remark}


\subsection{Truncation of functions}

We define the truncation operator and state its basic properties.

\begin{definition}\label{def:truncation}
For $M>0$, we define the truncation operator $\tr_M: L^1(\o;\R^M)\to L^\infty(\Omega;\R^M)$ as
\[
\tr_M(f)(x) \coloneqq \begin{dcases}
        f(x) & |f(x)|\leq M,\\
        M\frac{f(x)}{|f(x)|} & |f(x)|>M.
    \end{dcases}
\]
\end{definition}

The following can be easily proved.
\begin{lemma}\label{lem:truncation}
Let $f\in W^{1,2}(\Omega;\R^M)$. Then, $\tr_M(f) \in L^\infty(\Omega;\R^M) \cap W^{1,2}(\Omega;\R^M)$ and $|\nabla \tr_M(f)|\leq |\nabla f|$.
\end{lemma}


\subsection{Geodesics of degenerate metrics}
Here we describe results from \cite{SternbergZuniga}, that will be used extensively to prove convergence of the degenerate geodesic problems under mild assumptions.

Let $F:\R^M \to [0,\infty)$ be a continuous function satisfying
\begin{enumerate}[label= (F\arabic*)]
    \item \label{F1} The zero set of $F$, denoted by $\mathcal{Z}$, consists of a finite number of distinct points;
    \item \label{F2} $\liminf_{|z|\to\infty}F(z)>0$.
\end{enumerate}

As in \Cref{def:limiting_constant}, we define $\mathrm{Lip}_{\mathcal{Z}}([-1,1];\R^M)$ to be the space of continuous curves which are Lipschitz continuous with respect to the Euclidean metric on any compact portion of the curve that does not touch the zero set of $F$.

Consider the energy
\[
E(\gamma):=\int_{-1}^1 F(\gamma(t))|\gamma'(t)|dt, \quad\quad \mathrm{ for }\,\, \gamma \in \mathrm{Lip}_{\mathcal{Z}}([-1,1];\R^M).
\]
Due to the parameterization invariance of the energy and the fact that it is conformal to the Euclidean metric up to a degenerate factor, we can define a metric on $\R^M$ by 
\[
d(p,q) \coloneqq \inf\left\{\int_{-1}^1 {F(\gamma(t))}|\gamma'(t)|dt
    : \gamma\in \mathrm{Lip}_\mathcal{Z}([0,1];\R^M),\, \gamma(-1)=p,\, \gamma(1)=q \right\},
\]
and $(\R^M,d)$ is a length space. Further, we introduce the length functional $L$ for any curve $\gamma$ as
       \[
       L(\gamma) := \sup_{\{t_k\}_k \subset\mathcal{P}}{\sum_{k} d(\gamma(t_k),\gamma(t_{k+1}))},
       \]
       where $\mathcal{P}$ is the set of finite partitions of $[-1,1]$.
Using this length space viewpoint, \Cref{def:SZThms} below was proven in \cite{SternbergZuniga}.

\begin{proposition}\label{def:SZThms}\hfill
       \begin{enumerate}[label=(\arabic*), ref=\arabic*]
       \item \cite[Lemma 2.4]{SternbergZuniga}\label{1} Let B(x,r) denote the open ball centered at x and with radius r in the Euclidean metric. For every $\e>0$ such that $\mathcal{Z}\subset B(0,\frac{1}{\e})$, there is an $r_\e>0$ such that if $p,q \in B(0,\frac{1}{\e}) \cap \left(\bigcup_{z\in \mathcal{Z}}B(z,2\e)^c\right)$ then there is a $d$-minimizing curve, $\gamma^* \in \mathrm{Lip}_\mathcal{Z}([-1,1];\R^M)$, such that
       \[
       \gamma^*([-1,1])\cap \left(\bigcup_{z\in \mathcal{Z}}B(z,\e)\right) = \emptyset;
       \]
        \item \cite[Theorem 2.5]{SternbergZuniga}\label{2} For any $\gamma \in \mathrm{Lip}_\mathcal{Z}([-1,1];\R^M)$, we have
        \[
        L(\gamma) = E(\gamma);
        \]
        \item \cite[Theorem 2.6]{SternbergZuniga}\label{3} For every $p,q \in \R^M$, there is a minimizer $\gamma^* \in \mathrm{Lip}_\mathcal{Z}([-1,1];\R^M)$ which satisfies
           \[
           d(p,q)=E(\gamma^*) = L(\gamma^*);
           \]
        \item \cite[Proposition 2.7]{SternbergZuniga}\label{4} Given any partition $\{t_k\}$ of $[-1,1]$, a minimizer $\gamma^* \in \mathrm{Lip}_\mathcal{Z}([-1,1];\R^M)$ satisfies 
        \[
        L(\gamma^*) = \sum_k d(\gamma^*(t_k),\gamma^*(t_{k+1})).
        \]
       \end{enumerate}   
   \end{proposition}


\subsection{$\Gamma$-convergence}

In this section, we recall the definition and the basic properties of $\Gamma$-limits. Since in this paper we work in the setting of the metric space $L^1(\o;\R^M)$, we will present the equivalent definition with sequences.
We refer to \cite{Dalmasobook} (see also \cite{Braides}) for a complete study of $\Gamma$-convergence on topological spaces.

\begin{definition}\label{def:gc}
Let $(X,\mathrm{d})$ be a metric space, and let $\{F_n\}_n$ be a sequence of functionals $F_n:X\to[-\infty,+\infty]$. We say that $\{F_n\}_n$ $\Gamma$-converges to $F:X\to[-\infty,+\infty]$
with respect to the metric $\mathrm{d}$, if the followings hold:
\begin{itemize}
\item[(i)] ($\Gamma$-Liminf) For every $x\in X$ and every $\{x_n\}_n\subset X$ with $x_n\to x$, we have
\[
F(x)\leq\liminf_{n\to\infty} F_n(x_n),
\]
\item[(ii)] (Recovery sequence) For every $x\in X$, there exists $\{x_n\}_n\subset X$ such that
\[
\limsup_{n\to\infty} F_n(x_n)\leq F(x),
\]
and with $x_n\to x$.
\end{itemize} \vspace{0\baselineskip}
\end{definition}

\subsection{Sets of finite perimeter}

We recall the definition and some basic facts about sets of finite perimeter that are needed in the paper. For more details on the subject, we refer the reader to standard references, such as \cite{AFP, EG, Giusti, MaggiBook}.

\begin{definition}
Let $E\subset\R^N$ with $|E|<\infty$, and let $A\subset\R^N$ be an open set.
We say that $E$ has \emph{finite perimeter} in $A$ if
\[
P(E;A)\coloneqq\sup\left\{\, \int_E \mathrm{div}\varphi \,d x \,:\, \varphi\in C^1_c(A;\R^N)\,,\, \|\varphi\|_{L^\infty}\leq1  \,\right\}<\infty.
\]
\end{definition}

\begin{definition}
Let $a,b\in\R^M$. We define the space $BV(\Omega;\{a,b\})$ as the space of functions $u\in L^1(\Omega;\R^M)$ with $u(x)\in\{a,b\}$ for a.e. $x\in\Omega$, and such that the set $\{x\in\Omega : u(x)=a\}$ has finite perimeter in $\Omega$.
\end{definition}

\begin{definition}
Let $E\subset\R^N$ be a set of finite perimeter in the open set $A\subset\R^N$. We define $\partial^* E$, the \emph{reduced boundary} of $E$, as the set of points $x\in\R^N$ for which the limit
\[
\nu_E(x)\coloneqq -\lim_{r\to0}\frac{D\chi_E(B(x,r))}{|D\chi_E|(B(x,r))}
\]
exists and is such that $|\nu_E(x)|=1$.
The vector $\nu_E(x)$ is called the \emph{measure theoretic exterior normal} to $E$ at $x$.
\end{definition}

We recall part of the De Giorgi's structure theorem for sets of finite perimeter.

\begin{theorem}
Let $E\subset\R^N$ be a set of finite perimeter in the open set $A\subset\R^N$. Then,
\[
P(E, B) = \mathcal{H}^{N-1}(\partial^* E\cap B),
\]
for all Borel sets $B\subset A$.
\end{theorem}


\section{Technical results}

In this section we collect the main technical results that will be used in the proofs of the main theorems.


\subsection{Estimates for sequences with uniformly bounded energy} 

We start by finding bounds that will allow to compare the energy of a sequence $\{u_n\}_n$ with the energy of the unfolded sequence $\{T_{\delta_n u_n}\}_n$.

\begin{remark}
We first remark that any $\{u_n\}_n \subset W^{1,2}(\Omega;\R^M)$ with bounded energy, i.e.,
\[
\sup_{n\in\N} \mathcal{E}_n (u_n) \leq C,
\] satisfies the following energy estimate
\begin{equation}\label{eq:grad_est}
\|\nabla u_n\|^2_{L^2(\o;\R^{N\times M})}\leq \frac{C}{\e_n}.
\end{equation}
By the chain rule (see \Cref{rmk:unfolding}) and that the unfolding operator is a bounded operator, we can compute
$$\|\nabla_y \unf u_n\|^2_{L^2(\o;L^2(Q;\R^{N\times M}))} = \delta_n^2\| \unf (\nabla u_n)\|^2_{L^2(\o;L^2(Q;\R^{N\times M}))} \leq \delta_n^2 \|\nabla u_n\|^2_{L^2(\o;\R^{N\times M})}.$$
Thus, we have the useful estimate
\begin{equation}\label{eq:grad_est_two_scale}
\|\nabla_y \unf u_n\|^2_{L^2(\o;L^2(Q;\R^{N\times M}))} \leq C\frac{\delta_n^2}{\varepsilon_n}.
\end{equation}
\end{remark}

By slightly modifying the key Poincairé-type estimates in \cite{BradiesZeppieri}\cite{Hagerty} rewritten in terms of the unfolding operator we achieve the key estimate that will be used throughout the paper.

\begin{theorem}\label{thm:est_2}
Let $\{u_n\}_n \subset W^{1,2}(\Omega;\R^M)$ be such that
\[
\sup_{n\in\N} \mathcal{E}_n (u_n) \leq C.
\]
Then,
    \begin{equation} \label{eq:poincaire_estimate}
		\|\unf u_n - u_n\|^2_{L^2(\o;L^2(Q;\R^{M}))} \leq C\left(\|\nabla_y \unf u_n\|^2_{L^2(\o;L^2(Q;\R^{N\times M}))}+\int_{\o\setminus\ohn}|u_n-a|^2 dx\right).
		\end{equation}
Moreover, if $\sup_n\|u_n\|_\infty<+\infty$, then
        \begin{equation} \label{eq:key_new_estimate_wo_grad}
			\|\unf u_n - u_n\|_{L^2(\o;L^2(Q;\R^{M}))} \leq C\delta_n^{\frac{1}{2}}.
		\end{equation}	
In particular, this implies that
		\begin{equation} \label{eq:key_new_estimate}
			\|\unf u_n - u_n\|_{L^2(\o;L^2(Q;\R^{M}))}\|\nabla u_n\|_{L^2(\o;\R^{N\times M})} \leq C\left(\frac{\delta_n}{\e_n}\right)^{\frac{1}{2}}.
		\end{equation}		
\end{theorem}

\begin{proof}
\textbf{Step 1.} We first prove \eqref{eq:poincaire_estimate}.
For $x\in\Omega$, let
\[
(\unf u_n)_Q(x) \coloneqq \int_Q \unf u_n(x,y)\, dy.
\]
Using the triangle inequality, together with the inequality $(p+q)^2\leq 2(p^2+q^2)$, we get
\begin{equation}\label{eq:poincaire_estimate_1}
\|\unf u_n - u_n\|^2_{L^2} \leq 
    2\|\unf u_n - (\unf u_n)_Q\|^2_{L^2} + 2\| u_n - (\unf u_n)_Q\|^2_{L^2}.
\end{equation}
where the norm is the $L^2(\o;L^2(Q;\R^{M}))$ norm.
We estimate the latter term on the right-hand side of \eqref{eq:poincaire_estimate_1}.
We split the integral as
$$\int_\o{|u_n-(\unf u_n)_Q|^2 dx} = \int_\ohn |u_n-(\unf u_n)_Q|^2 dx + \int_{\o\setminus\ohn}|u_n-a|^2 dx.$$
 Using the unfolding operator similarly to \Cref{rmk:unfolding}, and since $\unf [(\unf u_n)_Q] = (\unf u_n)_Q$, we get
\begin{equation}\label{eq:poincaire_estimate_2}
\int_\o{|u_n-(\unf u_n)_Q|^2 dx} \leq \|\unf u_n - (\unf u_n)_Q\|^2_{L^2(\o;L^2(Q;\R^{M}))} + \int_{\o\setminus\ohn}|u_n-a|^2 dx.
\end{equation}
Now, we estimate the first term on the right-hand side of \eqref{eq:poincaire_estimate_1}.
By the Poincar\'{e}-Wirtinger inequality in the $y$-variable, for each $x\in\Omega$, we can estimate
$$
\int_Q|\unf u_n - (\unf u_n)_Q|^2\;dy \leq C\int_Q|\nabla_y \unf u_n|^2\;dy.
$$
Integrating over $\Omega$ we get the bound
\begin{equation}\label{eq:poincaire_estimate_3}
\|\unf u_n - (\unf u_n)_Q\|^2_{L^2(\o;L^2(Q;\R^{M}))} \leq C\|\nabla_y \unf u_n\|^2_{L^2(\o;L^2(Q;\R^{N\times M}))}.
\end{equation}
Thus, from \eqref{eq:poincaire_estimate_1}, \eqref{eq:poincaire_estimate_2}, and \eqref{eq:poincaire_estimate_3}, we deduce \eqref{eq:poincaire_estimate}.\\

\textbf{Step 2.} We now prove \eqref{eq:key_new_estimate_wo_grad} and \eqref{eq:key_new_estimate}.
Since, $\sup_n\|u_n\|_\infty<+\infty$, we get
$$\int_{\o\setminus\ohn}|u_n-a|^2 dx \leq C|\o\setminus\ohn|\leq C\delta_n.$$

Combining \eqref{eq:poincaire_estimate} together with \eqref{eq:grad_est_two_scale}, and that $\delta_n \ll \e_n$, we achieve \eqref{eq:key_new_estimate_wo_grad}
\[
\|\unf u_n - u_n\|_{L^2(\o;L^2(Q;\R^{M}))} \leq C\left(\frac{\delta_n^2}{\e_n}+ \delta_n \right)^{\frac{1}{2}} \leq C\delta_n^{\frac{1}{2}}.
\]

In view of \eqref{eq:grad_est}, we conclude that 
\[
\|\unf u_n - u_n\|_{L^2(\o;L^2(Q;\R^{M}))}\|\nabla u_n\|_{L^2(\o;\R^{N\times M})} \leq C\left(\frac{\delta_n}{\e_n}\right)^{\frac{1}{2}}.
\]

\end{proof}


\subsection{Definition and properties of the auxiliary cell problem}

In this section, we study an auxiliary cell problem that will be invoked in the proof of the liminf inequality (see Proposition \ref{thm:liminf}).

\begin{definition} \label{def:W_eta}Define the function $W^\eta:\R^M \to [0,\infty)$ as
 	$$W^\eta(z) := \inf_{\psi \in \mathcal{A}_\eta}\int_Q W(y,z + \psi(y))\;dy,$$ 
where the admissible set $\mathcal{A}_\eta$ is given by
\begin{align*}
&\mathcal{A}_\eta:= \Big\{\psi \in L^\infty(Q;\R^M) \cap W^{1,2}_0(Q;\R^M) : \|\psi\|_{L^\infty(Q;\R^M)}\leq \eta, \\
&\hspace{6cm}  \|\psi\|_{L^2(Q;\R^M)}\|\nabla \psi\|_{L^2(Q;\R^{N\times M})}\leq 5\eta^2 \Big\}.
\end{align*}
\end{definition}

We prove some properties of the function $W^\eta$.

\begin{theorem}[Properties of $W^\eta$]\label{thm:prop_W_eta}
The followings hold:
	\begin{enumerate}
		\item For every $z \in \R^M$, the infimum problem defining $W^\eta(z)$ admits a minimizer;
  		\item $W^\eta$ is continuous;
		\item $W^\eta(z) = 0 \iff z \in \{a,b\}$;
		\item For each $z\in\R^M$, $W^\eta(z)$ converges increasingly to \[
        W_{\mathrm{\mathrm{hom}}}(z) \coloneqq \int_Q W(y, z)\, dy,
        \]
        as $\eta \to 0$. Moreover, $W^\eta$ converges uniformly to $W_{\mathrm{\mathrm{hom}}}$ on every compact set.
	\end{enumerate}
\end{theorem}

\begin{proof}
 \textbf{Step 1.} We prove (1). Fix $z\in\R^M$. Let $\{\psi_n\}_n$ be an infimizing sequence for $W^\eta(z)$.
	Since $\sup_n \|\psi_n\|_\infty \leq \eta$, up to a subsequence (not relableled), we have that $\psi_n \wkstr \psi$ for some $\psi\in L^\infty(Q;\R^M)$, and, in turn, $\psi_n \wk \psi$ in $L^2(Q;\R^M)$.
    This is not enough to conclude by using the lower semicontinuity of the integral functional. We need to improve the convergence. To do that, we now consider two cases.\\
	\emph{Case 1.} Suppose $\psi\neq 0$. Using the constraints satisfied by each $\psi_n$'s, we get
 \[
\limsupn \|\nabla\psi_n\|_{L^2(Q;\R^M)} \leq \limsupn \frac{5\eta^2}{\|\psi_n\|_{L^2(Q;\R^M)}}
    \leq \frac{5\eta^2}{\|\psi\|_{L^2(Q;\R^M)}},
 \]
 where the last step is obtained by the fact that
 \[
\|\psi\|_{L^2(Q;\R^M)} \leq \liminfn \|\psi_n\|_{L^2(Q;\R^M)}.
 \]
	Thus, we deduce that $\{\psi_n\}_n$ is bounded in $W^{1,2}(Q;\R^M)$. By the Rellich–Kondrachov Theorem, we get that $\psi_n\to \psi$ strongly in $L^2(Q;\R^M)$, and weakly in $W^{1,p}(Q;\R^M)$ for all $p\in[1,2]$.
    In particular, by the compactness of the trace operator, we note that $\psi\in W^{1,2}_0(Q;\R^M)$. We conclude that $\psi \in \mathcal{A}_\eta$, because 
    \begin{align*}
        \|\psi\|_{L^2(Q;\R^M)}\|\nabla \psi\|_{L^2(Q;\R^M)}&\leq (\liminfn \|\psi_n\|_{L^2(Q;\R^M)})(\liminfn \|\nabla\psi_n\|_{L^2(Q;\R^M)})\\
        &\leq \liminfn \|\psi_n\|_{L^2(Q;\R^M)}\|\nabla\psi_n\|_{L^2(Q;\R^M)} \leq 5\eta^2.
    \end{align*}
	 
	\emph{Case 2.} Now we consider the case in which $\psi=0$.
    Note that if
    \[
    \liminfn \|\nabla\psi_n\|_{L^2(Q;\R^M)} < +\infty,
    \]
    we can argue as in the previous case, by taking a subsequence bounded in $W^{1,2}(Q;\R^M)$.
	 We now consider the case in which
  \[
  \liminfn \|\nabla\psi_n\|_{L^2(Q;\R^M)} = +\infty.
  \]
  We have
	 \[
  \limsupn \|\psi_n\|_{L^2(Q;\R^M)} \leq
  \limsupn \frac{5\eta^2}{\|\nabla\psi_n\|_{L^2(Q;\R^M)}}
  = 0.
  \]

  Thus, in both cases we achieve that $\psi_n\to \psi$ strongly in $L^2(Q;\R^M)$ and $\psi \in \mathcal{A}_\eta$. We then conclude by using the Dominated Convergence Theorem (thanks to the upper bound on $W$ \ref{W4}) and the uniform continuity of $W$ on $\overline{B(0,|z|+\eta)}$ that
  \[
 \lim_{n\to\infty}\int_Q W(y,z + \psi_n(y))\;dy = \int_Q W(y,z + \psi(y))\;dy.
  \]

\textbf{Step 2.} We establish (2). Let $\{z_n\}_n\subset\R^M$ be such that $z_n\to z$.
We first prove that
\[
W^\eta(z) \leq\liminfn W^\eta(z_n).
\]
Using Step 1, for each $n\in\N$ there exists $\psi_n\in \mathcal{A}_\eta$ such that
\[
W^\eta(z_n) = \int_Q W(y, z_n + \psi_n(y))\;dy,
\]
and clearly,
\begin{equation}\label{eq:W_conv_1}
W^\eta(z) \leq \int_Q W(y, z_n + \psi_n(y))\;dy
\end{equation}
for all $n\in\N$. Since $\| \psi_n \|_{L^\infty(Q;\R^M)}\leq \eta$ for all $n\in\N$, we have that
\begin{equation}\label{eq:W_conv_2}
\lim_{n\to\infty} \left| W(y, z_n + \psi_n(y)) - W(y, z + \psi_n(y)) \right| = 0,
\end{equation}
for all $y\in Q$. Using \eqref{eq:W_conv_1}, \eqref{eq:W_conv_2}, the upper bound on $W$ \ref{W4}, we can apply Dominated Convergence Theorem to conclude that
\[
W^\eta(z) \leq \liminfn \int_Q W(y, z_n + \psi_n(y))\;dy = \liminfn W^\eta(z_n),
\]
as desired.

In order to establish the other inequality, let $\psi \in \mathcal{A}_\eta$ be such that
\[
W^\eta(z) = \int_Q W(y, z + \psi(y))\;dy.
\]
Then, for each $n\in\N$, we get
\begin{align*}
W^\eta(z_n) &\leq \int_Q W(y, z_n + \psi(y))\;dy \\
& \leq \int_Q W(y, z + \psi(y))\;dy
    + \int_Q \left| W(y, z_n + \psi_n(y)) - W(y, z + \psi_n(y)) \right|\, dy \\
&\leq W^\eta(z)
    + \int_Q \left| W(y, z_n + \psi_n(y)) - W(y, z + \psi_n(y)) \right| \, dy.
\end{align*}
Thus, by \eqref{eq:W_conv_2}, we get the opposite inequality.\\

\textbf{Step 3.} We prove (3). This follows directly from the fact that $W\geq 0$, that (see \ref{W2}) $W(x,z)=0$ if and only if $z\in\{a,b\}$, and that $\psi=0$ on $\partial Q$ in the sense of traces.\\

\textbf{Step 4.} We prove (4). Let $\eta_1<\eta_2$. Then, $\mathcal{A}_{\eta_1}\subset \mathcal{A}_{\eta_2}$, and thus $W^{\eta_1}(z)\geq W^{\eta_2}(z)$ for all $z\in\R^M$. Let $\{\eta_n\}_n$ with $\eta_n\to0$ as $n\to\infty$. For each $n\in\N$, let $\psi_n\in \mathcal{A}_{\eta_n}$ be a solution for the minimizing problem defining $W^{\eta_n}(z)$.
We claim that $\psi_n\to 0$ strongly in $L^2$. Indeed, since $\|\psi_n\|_{L^\infty}(Q;\R^M)\leq \eta_n$, we get the claim.
Similar to the previous proofs, the upper bound on $W$ \ref{W4} allows us to apply Dominated Convergence Theorem to prove that $W^{\eta}(z)\to W_{\mathrm{\mathrm{hom}}}(z)$ as $\eta\to0$. Since the sequence is increasing, we can apply Dini's Theorem and get that the convergence is uniform on compact sets.
\end{proof}


\subsection{Properties of distances with degenerate metrics}

In this section we study the properties of the metrics
\[
d_\eta(p,q) \coloneqq \inf\left\{\int_{-1}^1 2\sqrt{W^\eta(\gamma(t))} |\gamma'(t)|dt
    : \gamma\in \mathrm{Lip}_\mathcal{Z}([-1,1];\R^M),\, \gamma(-1)=p,\, \gamma(1)=q \right\},
\]
where the space $\mathrm{Lip}_\mathcal{Z}([-1,1];\R^M)$ is introduced in Definition \ref{def:limiting_constant}.

In particular, we are interested in their behavior as $\eta\to0$.
We first state some properties that are easy to establish.

    \begin{lemma}\label{lem:limit_deg_metric}
        For $p,q \in \R^M$, define 
        \[
        d_0(p,q) := \sup_{\eta>0} d_\eta(p,q).
        \]
        Then
        \begin{enumerate}
            \item $d_0$ is a metric on $\R^M$;
            \item For each $p,q\in\R^M$, it holds $\lim_{\eta\to0} d_\eta(p,q) = d_0(p,q)$.
            \item $d_0 \leq d_{\mathrm{hom}}$ where
{$$d_{\mathrm{hom}}(p,q) \coloneqq \inf_{\gamma\in \mathrm{Lip}_\mathcal{Z}([-1,1];\R^M)}\left\{\int_{-1}^1 2\sqrt{W_{\mathrm{hom}}(\gamma(t))} |\gamma'(t)|dt
    : \gamma(-1)=p, \; \gamma(1)=q \right\}.$$}
\end{enumerate}
    \end{lemma}

\begin{proof}
\textbf{Step 1.} We first prove that $d_0$ is a metric on $\R^M$. Since each $d_\eta$ is a metric on $\R^M$, we get that $d_0(p,q)\geq0$ for all $p,q\in\R^m$, and that $d_0(p,q)=0$ if and only if $p=q$.
We now prove the triangle inequality. Let $p,q,r\in\R^M$. Then, invoking the triangle inequality for each $d_\eta$, we get
\begin{align*}
d_0(p,q) &= \sup_{\eta>0} d_\eta(p,q) \leq \sup_{\eta} [ d_\eta(p,r) + d_\eta(r,q) ] \\
&\leq \sup_{\eta>0} d_\eta(p,r) + \sup_{\eta} d_\eta(r,q)  = d_0(p,r) + d_0(r,q),
\end{align*}
as desired.

\textbf{Step 2.} For $0<\eta_1<\eta_2$, by Theorem \ref{thm:prop_W_eta} we have that $W^{\eta_1}(z)\geq W^{\eta_2}(z)$, for all $z\in\R^M$. Thus, for each $p,q\in\R^M$, the supremum in the definition of $d_0(p,q)$ is actually a limit.
\textbf{Step 3.} This follows directly from \Cref{thm:prop_W_eta} as $W^\eta \leq W_{hom}$.
\end{proof}
 
The existence of minimizing geodesics has already been established in \cite{SternbergZuniga} for each $\eta$ (see Proposition \ref{def:SZThms}(\ref{3})). Therefore, for each $\eta>0$ we have $\gamma_\eta  \in \mathrm{Lip}_\mathcal{Z}([0,1];\R^M)$ with $\gamma_\eta(-1)=a,\, \gamma_\eta(1)=b$ such that
\[
\sigma_\eta\coloneqq d_\eta(a,b) = \int_{-1}^1 2\sqrt{{W}^{\eta}(\gamma_\eta)}|\gamma_\eta'|dt.
\]
For notational convenience, we also define
$$\sigma_0\coloneqq \sup_{\eta>0}d_\eta(a,b).$$

We now  investigate the behavior of minimizing curves $\gamma_\eta$. 

    \begin{lemma}\label{lem_gamma_0}
        Let $\{\eta_n\}_n$ be an infinitesimal sequence, and for each $n\in\N$ let $\gamma_{\eta_n}$ be a geodesics for $d_{\eta_n}(a,b)$.
        Then, (up to a subsequence, not relabeled) there exists a curve $\gamma_0\in \mathrm{Lip}_\mathcal{Z}([-1,1];\R^M)$ such that
        \[
        \lim_{n\to\infty} \sup_{t\in[-1,1]} d_0(\gamma_{\eta_n}(t), \gamma_0(t))=0,
        \]
        and $\gamma_0(-1)=a$, $\gamma_0(1)=b$.
    \end{lemma}
    
\begin{proof}
We apply the Ascoli-Arzel\`{a} Theorem.

\textbf{Step 1.} We first prove equiboundedness with respect to the metric $d_0$. We claim that there exists $\widetilde{M}>0$ such that
\[
\{\gamma_{\eta_n}(t):\; t \in [-1,1]\} \subset B_{d_0}(a,\widetilde{M}),
\]
for all $n\in\N$.
We start by showing that there exists $M>0$ such that
\[
\{\gamma_{\eta_n}(t):\; t \in [-1,1]\} \subset B(0,M),
\]
for all $n\in\N$, where $B(0,M)$ is the Euclidean ball.
Indeed, note that
\begin{equation}\label{eq:bound_geodesics_M}
d_{\eta_n}(a,b) \leq d_{hom}(a,b),
\end{equation}
independent of $n$.
Since $W$ grows linearly at infinity (see \ref{W3}), we get that there exists $\widetilde{R}>0$ such that
\begin{equation}\label{eq:est_W_quad}
W^{\eta_n}(z)\geq C(d_{\mathrm{hom}}(a,b))^2|z|,
\end{equation}
for all $|z|\geq \widetilde{R}$, where the constant $C>0$ is independent of $n$ and of $z\in\R^M$.
Let $M> \max\{|a|, |b|, \widetilde{R}\}$.
Assume that there exists $t_n\in[-1,1]$ such that $|\gamma_n(t_n)|>2M$.
Since the curve $\gamma_n$ is continuous and $\gamma_n(-1), \gamma_n(1)\in B(0,M)$, by the choice of $M$, it is possible to find $t^n_1(M) < t_n < t^n_2(M)$ such that
$|\gamma_n(t)|\geq M$ for all $t\in [t_1^n(M), t_2^n(M)]$, and
$|\gamma_n(t^n_1(M))| = |\gamma_n(t^n_2(M))| = M$.
Then, by using \eqref{eq:est_W_quad}, we get
\begin{align*}
d_{\eta_n}(a,b) & = 2\int_{-1}^1 \sqrt{W(\gamma_n(t))} |\gamma'_n(t)| dt \\
&\geq C d_{hom}(a,b)\sqrt{M} \int_{t^n_1(M)}^{t^n_2(M)} |\gamma'_n(t)| dt
\\
&\geq C d_{hom}(a,b) M^{3/2},
\end{align*}
where in the third step we used the choice of $t^n_1(M)$ and $t^n_2(M)$.
This contradicts \eqref{eq:bound_geodesics_M} for $M$ large enough.

We now conclude this step as follows. Using the upper bound of $d_0$ by $d_\mathrm{hom}$ and choosing the straight line path between the point $a$ and the termination point (denoted by $L_{a,p}$), we can compute
\[
\sup_{t \in [-1,1]}d_0(a,\gamma_{\eta_n}(t))\leq \sup_{t \in [-1,1]}d_\mathrm{hom}(a,\gamma_{\eta_n}(t))\leq \sup_{t \in [-1,1]}2|\gamma_{\eta_n}(t)-a|\int_{-1}^{1}\sqrt{W_{\mathrm{hom}}(L_{a,\gamma_{\eta_n}(t)}})ds .
\]
using the continuity of $\sqrt{W_{\mathrm{hom}}}$ and the boundedness of $\gamma_{\eta_n}(t)$ in the Euclidean ball $B(0,M)$, we can establish the uniform boundedness in the $d_0$ metric
\[
\{\gamma_{\eta_n}(t):\; t \in [0,1]\} \subset \overline{B_{d_0}(a,4\|\sqrt{W_{\mathrm{hom}}}\|_{L^\infty(B(0,M))}(|M|+|a|))}.
\]
We conclude by setting $\widetilde{M}\coloneqq 1+4\|\sqrt{W_{\mathrm{hom}}}\|_{L^\infty(B(0,M))}(|M|+|a|)$.

\textbf{Step 2.} We now prove the equicontinuity of the sequence with respect to the metric $d_0$. 
As in \cite[Proof of Theroem 2.6]{SternbergZuniga}, for any $n\in\N$ we rescale the curve $\gamma_\eta$ to a curve (that with an abuse of notation we still denote by) $\gamma_{\eta_n}:[-1, 1]$ in such a way that
\[
2\sqrt{W^{\eta_n}(\gamma_{\eta_n}(t))}|\gamma_{\eta_n}'(t)|\equiv \frac{\sigma_{\eta_n}}{2},
\]
for all $t\in [-1, 1]$.
Let $t_2>t_1$. Then,
\begin{align*}
d_{\eta_n}(\gamma_{\eta_n}(t_1),\gamma_{\eta_n}(t_2)) =
    \int_{t_1}^{t_2} 2\sqrt{W^{\eta_n}(\gamma_{\eta_n}(t))} |\gamma'_{\eta_n}(t)|\, dt
    =  \frac{\sigma_{\eta_n}}{2} (t_2-t_1)
    \leq \frac{\sigma_0}{2} (t_2-t_1).
\end{align*}
Therefore, we get that
\[
d_0(\gamma_{\eta_n}(t_1),\gamma_{\eta_n}(t_2))\leq \frac{\sigma_0}{2} (t_2-t_1) +\|d_{\eta_n}-d_0\|_\infty,
\]
where the last norm is the uniform norm in the space $\overline{B(0,M)}\times\overline{B(0,M)}$. Note that by definition of $d_0$, and by using Dini's Theorem on the compact space $\overline{B(0,M)}\times\overline{B(0,M)}$, we get that
\[
\lim_{n\to\infty} \|d_{\eta_n}-d_0\|_\infty=0.
\]
Fix $\varepsilon>0$, choose $\bar{n}\in\N$ such that
\[
\|d_{\eta_n}-d_0\|_\infty < \frac{\varepsilon}{2},
\]
and let
\[
\delta_0 \coloneqq \frac{\varepsilon}{\sigma_0}.
\]
Then, for every $n\geq \bar{n}$, it holds
\[
d_0(\gamma_{\eta_n}(t_1),\gamma_{\eta_n}(t_2))\leq \varepsilon,
\]
whenever $|t_2-t_1|<\delta_0$. For each $n=1,\dots,\bar{n}$, using the equiboundedness and uniform continuity of $\gamma_{\eta_n}$, let $\delta_n>0$ be such that
\[
d_0(\gamma_{\eta_n}(t_1),\gamma_{\eta_n}(t_2))\leq \varepsilon
\]
whenever $|t_2-t_1|<\delta_n$.
Define
\[
\delta\coloneqq \min\{ \delta_0, \delta_1,\dots,\delta_{\bar{n}} \}.
\]
This proves that the sequence $\{\gamma_{\eta_n}\}_{n\in\N}$ is equicontinuous.
\\

\textbf{Step 3.}
Applying the Ascoli-Arzel\`{a} Theorem, we get the existence of a subsequence and of a curve $\gamma_0:[-1,1]\to\R^M$ to which the subsequence converges uniformly with respect to the metric $d_0$.\\

\textbf{Step 4.}
We now prove that $\gamma_0\in \mathrm{Lip}_\mathcal{Z}([-1,1];\R^M)$.
This follows directly from the fact that, given any $r>0$, there exists $m_r>0$ such that
\[
\inf\left\{\, W^{\eta_n}(z) \,:\, z\in\R^M\setminus \left( B(a,r)\cup B(b,r) \right) \,\right\}\geq m_r,
\]
for all $n\in\N$.
Indeed, this can be deduced by using the uniform convergence of $W^{\eta_n}$ to $W_{\mathrm{hom}}$, together with the fact that $W_{\mathrm{hom}}$ only vanishes at $a$ and $b$.
\end{proof}

Now we are ready to prove the main result of this section.

   \begin{proposition}\label{prop:d0_character}
   Let $\gamma_0\in \mathrm{Lip}_\mathcal{Z}([-1,1];\R^M)$ be the curve given by Lemma \ref{lem_gamma_0}.
   Then, it holds
       \[
        d_0(a,b) =  \lim_{\eta \to 0} \int_{-1}^1 {2\sqrt{W^\eta(\gamma_0)}|\gamma_0'|dt}
            = \sup_{\eta>0} \int_{-1}^1 {2\sqrt{W^\eta(\gamma_0)}|\gamma_0'|dt}.
        \]
   \end{proposition}

   \begin{proof}
     By definition, we have
     \[
     d_\eta(a,b)\leq \int_{-1}^1 {2\sqrt{W^\eta(\gamma_0)}|\gamma_0'|dt},
     \]
     thus taking the limit as $\eta\to0$, and recalling that $\eta\mapsto W^\eta(z)$ is monotone for each $z\in\R^M$, we get
     \[
        d_0(a,b) \leq \lim_{\eta \to 0} \int_{-1}^1 {2\sqrt{W^\eta(\gamma_0)}|\gamma_0'|dt}.
     \]
     Now we prove the converse inequality. Using Proposition \ref{def:SZThms} (\ref{2}) and the fact that $\gamma_0\in \mathrm{Lip}_\mathcal{Z}([-1,1];\R^M)$, we get
     \[
     \int_{-1}^1 {2\sqrt{W^\eta(\gamma_0)}|\gamma_0'|dt} = L_\eta(\gamma_0).
     \]
     Fix an arbitrary finite partition $\{t_k\}_{k=1}^m$ of $[-1.1]$.
     Use the triangle inequality and the definition of $d_0$ to get
     \begin{align*}
         d_\eta(\gamma_0(t_k),\gamma_0(t_{k+1})) &\leq d_\eta(\gamma_\eta(t_k),\gamma_0(t_{k}))+d_\eta(\gamma_\eta(t_k),\gamma_\eta(t_{k+1}))+d_\eta(\gamma_0(t_k),\gamma_\eta(t_{k+1}))\\
         &\leq d_0(\gamma_\eta(t_k),\gamma_0(t_{k}))+d_\eta(\gamma_\eta(t_k),\gamma_\eta(t_{k+1}))+d_0(\gamma_0(t_k),\gamma_\eta(t_{k+1})).
     \end{align*}
        Fix $j \in \N$. In view of the uniform convergence given by \Cref{lem_gamma_0}, we can find $\eta_0(j,m)$ such that for all $\eta<\eta_0$ we have
        \[
        \sup_{t\in[-1,1]} d_0(\gamma_\eta(t),\gamma_0(t))\leq \frac{1}{jm}.
        \]
     Thus, we can bound the total sum over the partition using the uniform convergence estimate for all $\eta<\eta_0$ by
     \[
     {\sum_{k=1}^m d_\eta(\gamma_0(t_k),\gamma_0(t_{k+1}))} \leq
        \sum_{k=1}^m d_\eta(\gamma_\eta(t_k),\gamma_\eta(t_{k+1})) + \frac{1}{j} = L_\eta(\gamma_\eta) + \frac{1}{j} = d_\eta(a,b) + \frac{1}{j}.
        \]
     Taking a supremum over all possible finite partitions, we have
     \[
     \int_{-1}^1 {2\sqrt{W^\eta(\gamma_0)}|\gamma_0'|dt} = L_\eta(\gamma_0) \leq d_\eta(a,b) + \frac{1}{j},
     \]
     and we conclude taking the limit as $\eta \to 0$ and then $j \to \infty$.
   \end{proof}


\section{Liminf Inequality}

The goal of this section is to prove the following result.

\begin{theorem}\label{thm:liminf}
Let $\{u_n\}_n \subset W^{1,2}(\Omega;\R^M)$ be such that $u_n \to u\in BV(\o;\{a,b\})$ strongly in $L^1(\Omega;\R^M)$.
Then,
\[
\liminfn \mathcal{E}_n(u_n)\geq  \mathcal{E}_\infty(u).
\]
\end{theorem}

Note that this proposition is weaker than the liminf inequality, since we are assuming the limiting function to be in $BV(\o ;\{a,b\})$, and not just in $L^1(\Omega;\{a,b\})$. We will prove in the next section that this is sufficient. The reason why we are first proving the above result is that, in the proof of the compactness (see Theorem \ref{thm:main_compactness}), we will use arguments that are the core of the idea to prove the above result. Thus, we opt by presenting them first here.

Without loss of generality, in what follows we assume that
\[
\liminfn \mathcal{E}_n(u_n) = \lim_{n\to\infty} \mathcal{E}_n(u_n) < \infty.
\]

\textbf{Step 1: Reduction} In this step, we show that it suffices to prove \Cref{thm:liminf} for $u_n$ uniformly bounded in $L^\infty$, and $W$ being linear outside a ball.
Let $M>R$, where $R>0$ is the constant given in \ref{W3}.
Let $\psi_M:[0,+\infty)\to[0,1]$ be a smooth function such that
\[
\psi_M\equiv 1 \text{ on } [0,M],\quad\quad\quad
\psi_M\equiv 0 \text{ for } t\geq 2M.
\]
Define
\[
\widetilde{W}^M(x,z)\coloneqq \psi_M(|z|)W(x,z) + (1-\psi_M(|z|))\frac{|z|}{R},
\]
and
\[
\widetilde{\mathcal{E}}_n(v)\coloneqq \int_\o{\left[\frac{1}{\e_n}\widetilde{W}^M \left(\frac{x}{\delta_n}, v(x)\right) + \e_n|\nabla v(x)|^2 \right]\;dx}.
\]
We claim the following. Assume that for all $v\in BV(\Omega;\{a,b\})$ it holds
\begin{equation}\label{eq:ineq_vn}
\mathcal{E}_\infty(v)\leq \liminfn \widetilde{\mathcal{E}}_n(v_n),
\end{equation}
whenever $\{v_n\}_n\subset W^{1,2}(\Omega;\R^M)$ with $\|v_n\|_{L^\infty}\leq 2M$ is such that $v_n\to v$ strongly in $L^1(\o;\R^M)$. Then, \Cref{thm:liminf} holds.

Indeed, let $\{u_n\}_n \subset W^{1,2}(\Omega;\R^M)$ be such that $u_n \to u\in BV(\o;\{a,b\})$ strongly in $L^1(\o;\R^M)$. Let
\[
v_n\coloneqq \tr_{2M} u_n,
\]
where $\tr_{2M} u_n$ is the truncation of $v$ (see Definition \ref{def:truncation}).
We observe that
\begin{equation}\label{eq:conv_trunc}
\lim_{n\to\infty}\|u_n-v_n\|_{L^1(\o;\R^M)} = 0.
\end{equation}
Indeed, let
\[
\Delta_n:=\{x \in \o\;:\; |u_n(x)|>2M\}.
\]
By Chebyshev's inequality, we get
\begin{align}\label{eq:trunc_un}
    \mathcal{L}^N(\Delta_n) \leq \frac{1}{2M}\int_{\Delta_n}{|u_n|\;dx}
    \leq \frac{R}{2M}\int_{\o}{W\left(\frac{x}{\delta_n},u_n\right)\;dx}
    \leq C\frac{R}{2M}\e_n,
\end{align}
where in the second inequality we used the fact that $M>R_1$ together with \ref{W3} and in the last inequality we are using that the sequence has bounded energy.
Invoking the triangle inequality and convexity, we obtain by \eqref{eq:trunc_un}
\begin{align*}
    \|u_n-v_n\|_{L^1(\o;\R^M)} &= \int_{\Delta_n}|u_n-v_n|\;dx \\
&\leq \left(\int_{\Delta_n}{|u_n|\;dx} + 2M\mathcal{L}^N(\Delta_n)\right) \\
&\leq C_M\e_n,
\end{align*}
This proves \eqref{eq:conv_trunc}.

Therefore, by \eqref{eq:ineq_vn} we get
\[
\mathcal{E}_\infty(u)\leq \liminfn \widetilde{\mathcal{E}}_n(v_n),
\]
and to conclude, we claim that for all $n\in\N$
\[
\widetilde{\mathcal{E}}_n(v_n) \leq \mathcal{E}_n(u_n).
\]
Indeed, note that $\widetilde{W}_M$ satisfies \ref{W1}, \ref{W2}, and \ref{W3}. Thus, \ref{W3}, the definition of $v_n$, and the definition of $\widetilde{W}_M$ yield
\[
\widetilde{W}^M(x,v_n(x)) \leq \widetilde{W}^M(x,u_n(x)) \leq W(x, u_n),
\]
for almost all $x\in\Omega$, and all $z\in\R^M$.

This, together with Lemma \ref{lem:truncation}, yields
\[
\int_\o{\left[\frac{1}{\e_n}\widetilde{W}^M \left(\frac{x}{\delta_n},v_n(x)\right) + \e_n|\nabla v_n(x)|^2 \right]\;dx}\leq \mathcal{E}_n(u_n),
\]
and thus the desired conclusion.\\

In the rest of this section, we will therefore assume that
\begin{equation}\label{eq:assumptions_un}
\sup_{n\in\N} \|u_n\|_{L^\infty(\Omega;\R^M)} \leq M,
\end{equation}
for some $M>R$.


\textbf{Step 2: Slicing}
For each $n\in\N$ and $\eta>0$, let
\[
f^\eta_n(x,y)\coloneqq \tr_\eta(\unf u_n(x,y) - u_n(x)),
\]
where the truncation operator $\tr_\eta$ is with respect to the variable $y$.
We want use a De Giorgi's slicing type of argument to modify $f^\eta_n$ to make it vanish on $\Omega\times\partial Q$.

For $k\in \N\setminus\{0\}$, let (see Figure \ref{fig:slicing})
\[
S_{k} \coloneqq \left\{y\in Q \;:\;\frac{1}{k} <\text{dist}(y,\partial Q)\leq \frac{2}{k} \right\},
\]
and
\[
Q_{k} \coloneqq \left\{y\in Q \;:\; \frac{2}{k} <\text{dist}(y,\partial Q) \right\}.
\]
Let $\psi_k:Q\to[0,1]$ be a smooth function with $0\leq \psi\leq 1$ such that
\[
\psi_k\equiv 1 \text{ in } Q_{k},\quad\quad
\psi_k\equiv 0 \text{ on } Q\setminus (Q_{k}\cup S_k),\quad\quad
\|\nabla\psi\|_{L^\infty(Q)} \leq k.
\]
Define
\[
v^\eta_{n,k}(x,y) = \psi_{k}(y) f^\eta_n(x,y).
\]

\begin{figure}
\includegraphics[scale=0.7]{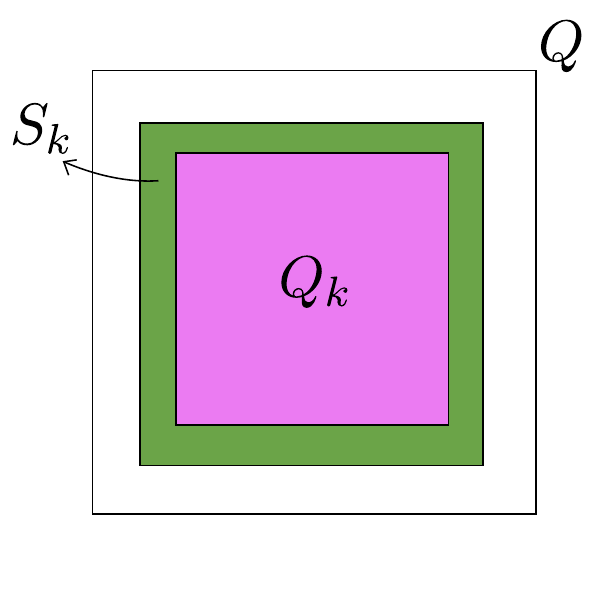}
\caption{The construction of the sets $Q_k$, and $S_k$, where the cut-off function $\psi_k$ equals one, and transition from one to zero, respectively.}
\label{fig:slicing}
\end{figure}

\begin{theorem}\label{thm:sequence_kxn}
There exists a sequence $\{k_n^\eta(x)\}_n$ with $x\mapsto k_n^\eta(x)$ measurable,
such that
\[
\limn k_n^\eta(x) = +\infty
\]
for almost every $x \in \o$,
\begin{equation}\label{eq:condition_kn_en}
\limsupn \int_\o{\frac{1}{k_n(x) \e_n}\;dx} = 0,
\end{equation}
and with $v_n^\eta\coloneqq v^\eta_{n,k_n(x)}$, it holds
	\begin{equation}
		\limsupn \mathcal{E}_n(u_n)\geq \limsupn\int_\o{2\left[\int_Q W(y,u_n+v_n^\eta)\;dy\right]^{\frac{1}{2}}|\nabla u_n|\;dx}.
	\end{equation}
Moreover, for almost every $x \in \o$, we have the estimate
    \begin{equation}\label{eqn:constraint_gradient}
		\|v_n^\eta(x,\cdot)\|_{L^2(Q;\R^M)}\|\nabla_y v_n^\eta(x,\cdot)\|_{L^2(Q;\R^M)}\leq
    4\eta^2 + \eta\| \nabla_y \unf u_n(x,\cdot) \|_{L^2(Q;\R^{N\times M})}.
	\end{equation}
\end{theorem}

The proof of the above theorem requires some preliminary results.
First, we claim that, for almost all $x\in\Omega$, it holds that
\begin{equation}\label{eq:est_v_eta_n_k}
\|v^\eta_{n,k}(x,\cdot)\|_{L^\infty(Q;\R^M)}\leq \eta,
\end{equation}
for all $n \in \N$, and that
\begin{equation}\label{eq:conv_v_eta_n_k}
\lim_{n\to\infty}\|v^\eta_{n,k}(x,\cdot)\|_{L^2(Q;\R^M)} 
=\lim_{n\to\infty}\|f^\eta_{n}(x,\cdot)\|_{L^2(Q;\R^M)}=  0.
\end{equation}
Indeed, the first estimate is direct from the definition of the function $v^\eta_{n,k}$ together with the fact that the cut-off function is bounded above by one. To see \eqref{eq:conv_v_eta_n_k}, we  observe that \eqref{eq:key_new_estimate_wo_grad} and the definition of $f^\eta_n$ gives that 
$$\lim_{n\to\infty}\|f^\eta_{n}\|_{L^2(\o;L^2(Q;\R^{M}))} = \lim_{n\to\infty}\|\unf u_n - u_n\|^2_{L^2(\o;L^2(Q;\R^{M}))} = 0.$$
In particular, this implies by taking a subsequence without relabeling that the following pointwise convergence holds for $\mathcal{L}^N a.e. \;x \in \o$ 
$$\lim_{n\to\infty}\|f^\eta_{n}(x,\cdot)\|_{L^2(Q;\R^M)}= 0.$$
Thus, \eqref{eq:conv_v_eta_n_k} follows from the fact that the cut-off function is bounded above by one. 

We can compute
		\begin{align}\label{eq:estimate_energy}
			&\int_\o{2\left[\int_Q W(y,u_n+ v^\eta_{n,k})\;dy\right]^{\frac{1}{2}}|\nabla u_n|\;dx} \nonumber \\
    &\hspace{2cm}\leq \int_\o{2\left[\int_Q W(y,u_n+f^\eta_n)\;dy\right]^{\frac{1}{2}}|\nabla u_n|\;dx} \nonumber \\
			&\hspace{3cm}\quad + \int_\o{2\left[\int_{Q\setminus Q_{k(x)}} W(y,u_n)\;dy\right]^{\frac{1}{2}}|\nabla u_n|\;dx} \nonumber \\
			&\hspace{3cm}\quad + \int_\o{2\left[\int_{S_{k(x)}} W(y,u_n+v^\eta_{n,k}(x,y))\;dy\right]^{\frac{1}{2}}|\nabla u_n|\;dx} \nonumber \\
			&\hspace{2cm}=: I_n + II_n + III_n.
	\end{align}
We will estimate the three integrals separately.

\begin{lemma}[Estimate for ${I}_n$] \label{lem:I_estimate}
It holds
\begin{align*}
   \limsupn &\int_\o 2 \left[\int_Q W(y,u_n+f^\eta_n)\;dy\right]^{\frac{1}{2}}|\nabla u_n|\;dx\nonumber\\
   &\leq \limsupn \int_\o 2 \left[\int_Q W(y,\unf u_n)\;dy\right]^{\frac{1}{2}}|\nabla u_n|\;dx.
\end{align*}
\end{lemma}

\begin{proof}
Write
\begin{align}\label{eq:writing_first_lemma}
&\int_\o 2\left[\int_Q W(y,u_n+\tr_\eta(\unf u_n - u_n))\;dy\right]^{\frac{1}{2}}|\nabla u_n|\;dx \nonumber \\
&\hspace{2cm}\leq \int_\o 2 \left[\int_Q W(y,\unf u_n)\;dy\right]^{\frac{1}{2}}|\nabla u_n|\;dx \nonumber\\
&\hspace{3cm}+ \int_\o 2 \left[\int_{Q_n^\eta(x)} W(y,u_n+\tr_\eta(\unf u_n - u_n))\;dy\right]^{\frac{1}{2}}|\nabla u_n|\;dx,
\end{align}
where
\[
Q_n^\eta(x):= \{y \in Q \;:\; |\unf u_n(x,y) - u_n(x)|>\eta\}.
\]
By Chebyshev inequality, we have
\begin{equation}\label{eq:second_writing_first_lemma}
\mathcal{L}^N(Q_n^\eta(x)) \leq \frac{1}{\eta^2}\int_Q{|\unf u_n(x,y) - u_n(x)|^2\;dy} = \frac{1}{\eta^2} \|\unf u_n(x,\cdot) - u_n(x)\|^2_{L^2(Q;\R^M)}.
\end{equation}
Using \ref{W4}, \eqref{eq:assumptions_un}, and \eqref{eq:second_writing_first_lemma} we get
\begin{align*}
    &\int_\o\left[\int_{Q_n^\eta(x)} W(y,u_n+\tr_\eta(\unf u_n - u_n))\;dy\right]^{\frac{1}{2}}|\nabla u_n|\;dx \\
    &\hspace{2cm}\leq \frac{C_{M}}{\eta^2}\int_\o{\mathcal{L}^N(Q_n^\eta(x))^{\frac{1}{2}}|\nabla u_n(x)|\;dx}\\
    &\hspace{2cm}\leq  \frac{C_{M}}{\eta^2}\|\mathcal{L}^N(Q_n^\eta(x))^{\frac{1}{2}}\|_{L^2(\o\;\R^M)}\|\nabla u_n\|_{L^2(\o\;\R^{N\times M})} \\
    &\hspace{2cm}= \frac{C_{M}}{\eta^2}\|\unf u_n - u_n\|_{L^2(\o;L^2(Q;\R^{M}))}\|\nabla u_n\|_{L^2(\o;\R^{N\times M})} \\
    &\hspace{2cm}\leq \frac{C_{M}}{\eta^2}
    \left(\frac{\delta_n}{\e_n}\right)^{\frac{1}{2}},
\end{align*}
where we used H\"{o}lder and  estimate \eqref{eq:key_new_estimate}. 
Since $\delta_n\ll\varepsilon_n$, in view of \eqref{eq:second_writing_first_lemma}this gives
\begin{align*}
\limsupn &\int_\o 2 \left[\int_Q W(y,u_n+\tr_\eta(\unf u_n - u_n))\;dy\right]^{\frac{1}{2}}|\nabla u_n|\;dx\nonumber\\
   &\leq \limsupn \int_\o 2 \left[\int_Q W(y,\unf u_n)\;dy\right]^{\frac{1}{2}}|\nabla u_n|\;dx.
\end{align*}
\end{proof}

Now, we estimate $II_n$.

\begin{lemma}[Estimate for ${II}_n$] \label{lem:II_estimate}
Assume that
\begin{equation}\label{eq:assumption_lemma_2}
\limsupn \int_\o{\frac{1}{k_n(x) \e_n}\;dx} = 0.
\end{equation}
Then,
\[
\lim_{n\to\infty} II_n = 0.
\]
\end{lemma}

\begin{proof}
Using \ref{W4} and the fact that $\|u_n\|_\infty\leq M$ (recall \eqref{eq:assumptions_un}) for all $n\in\N$, we get
	\begin{align*}
		\int_\o{2\left[\int_{Q\setminus Q_{k_n(x)}} W(y,u_n)\;dy\right]^{\frac{1}{2}}|\nabla u_n|\;dx} &\leq C_M \int_\o{\mathcal{L}^N(Q\setminus Q_{k_n(x)})^{\frac{1}{2}}|\nabla u_n|\;dx}\\
		&\leq C_M \|\nabla u_n\|_{L^2(\o)}\,
  \left(\int_\Omega\mathcal{L}^N(Q\setminus Q_{k_n(x)})\, dx\right)^{\frac{1}{2}},
	\end{align*}
where the last step follows from H\"{o}lder inequality.
Using estimate \eqref{eq:grad_est} together with the fact that
\[
\mathcal{L}^N(Q\setminus Q_{k_n(x)}) \leq  \frac{C}{k_n(x)},
\]
	we get
	$$II_n \leq C_M\left[\int_\o{\frac{1}{k_n(x) \e_n}\;dx}\right]^{\frac{1}{2}}. $$
Thanks to \eqref{eq:assumption_lemma_2}, we conclude.
\end{proof}

We finally estimate $III_n$ by employing a similar strategy.

\begin{lemma}[Estimate for ${III}_n$] \label{lem:III_estimate}
Assume that
\begin{equation}\label{eq:assumption_lemma_3}
\limsupn \int_\o{\frac{1}{k_n(x) \e_n}\;dx} = 0.
\end{equation}
Then,
\[
\lim_{n\to\infty} III_n = 0.
\]
\end{lemma}

\begin{proof}
Using \ref{W4} and the fact that $\|u_n+v^\eta_{n,k}\|_\infty\leq M+\eta$ for all $n\in\N$ (see \eqref{eq:assumptions_un}\eqref{eq:est_v_eta_n_k}), we get
\begin{align*}
&\int_\o{2\left[\int_{S_{k_n(x)}} W(y,u_n+v^\eta_{n,k}(x,y))\;dy\right]^{\frac{1}{2}}|\nabla u_n|\;dx}
\leq  C_{M,\eta}\int_\o{
 \mathcal{L}^N(S_{k_n(x)})^{\frac{1}{2}} |\nabla u_n|\;dx}.
\end{align*}

Applying H\"{o}lder inequality, we get
\[
III_n\leq C_{M,\eta} \|\nabla u_n\|_{L^2(\o)}
\left(\int_\Omega\mathcal{L}^N(S_{k_n(x)})\, dx\right)^{\frac{1}{2}}.
\]
and since
\[
\mathcal{L}^N(S_{k(x)}) \leq  \frac{C}{k_n(x)},
\]
by \eqref{eq:assumption_lemma_3} and \eqref{eq:grad_est} we conclude.
\end{proof}

We are now in position to prove the main theorem of this step.

\begin{proof}[Proof of Theorem \ref{thm:sequence_kxn}]
\underline{\emph{Part 1.}} Define
\[
k_n(x) \coloneqq \left\lfloor \frac{\eta^2}{\|f^\eta_n(x,\cdot)\|^2_{L^2(Q;\R^M)}}\right\rfloor .
\]
Observe that $k_n\geq 1$ is an integer for all $x\in\Omega$, and that the function $x\mapsto k_n(x)$ is Lebesgue measurable, as it is the composition of an upper semicontinuous function (hence Borel measurable) and a Lebesgue measurable function.
Furthermore, by \eqref{eq:conv_v_eta_n_k} we have that the denominator converges to zero, and thus $k_n \to \infty$ as $n\to\infty$.
Moreover, in view of the definition of $f^\eta_n$ and \eqref{eq:key_new_estimate_wo_grad}, we have
\[
\int_\o{\frac{1}{\e_n}\|f^\eta_n(x,\cdot)\|^2_{L^2(Q;\R^M)}\;dx}\leq \frac{1}{\e_n}\|\unf u_n - u_n\|^2_{L^2(\o;L^2(Q;\R^{M}))} \leq C\frac{\delta_n}{\e_n}.
\]
Thus, using that $\delta_n \ll \e_n$, we deduce that 
\[
\limsupn \int_\o{\frac{1}{k_n(x) \e_n}\;dx} = 0.
\]

\underline{\emph{Part 2.}} We now prove the energy estimate . From \eqref{eq:estimate_energy}, together with Lemmas \ref{lem:I_estimate}, \ref{lem:II_estimate}, and \ref{lem:III_estimate} we get the desired estimate
\begin{align*}
		&\limsupn\int_\o{2\left[\int_Q W(y,u_n+v_n^\eta)\;dy\right]^{\frac{1}{2}}|\nabla u_n|\;dx}\\
  &\hspace{2cm}\leq \limsupn \int_\o 2 \left[\int_Q W(y,\unf u_n)\;dy\right]^{\frac{1}{2}}|\nabla u_n|\;dx \\
  &\hspace{2cm}\leq \limsupn \mathcal{E}_n(u_n)
\end{align*}

\underline{\emph{Part 3.}} Finally, we establish the bound \eqref{eqn:constraint_gradient}.
Note that
\begin{equation*}
	|\nabla_y v_{n}^\eta|^2
    = |\nabla_y \psi_{k_n} f^\eta_n + \psi_{k_n(x)}\nabla_y f^\eta_n|^2
    \leq 2(k_n^2(x)|f^\eta_n|^2+|\nabla_y f^\eta_n|^2),
\end{equation*}
and thus
\begin{equation*}
	\|\nabla_y v_{n}^\eta (x,\cdot)\|_{L^2(Q;\R^M)}
    \leq \sqrt{2} \left[ k_n(x)\|f^\eta_n(x,\cdot)\|_{L^2(Q;\R^M)}
    +\|\nabla_y f^\eta_n(x,\cdot)\|_{L^2(Q;\R^M)} \right].
\end{equation*}
Since
\[
\|v_{n}^\eta\|_{L^2(Q;\R^M)} \leq \|f^\eta_n\|_{L^2(Q;\R^M)},
\]
we obtain
\begin{align*}
&\|v_{n}^\eta(x,\cdot)\|_{L^2(Q;\R^M)}\|\nabla_y v_{n}^\eta (x,\cdot)\|_{L^2(Q;\R^M)}\\
&\hspace{2cm}\leq
	 \sqrt{2} \left[ k_n(x)\|f^\eta_n(x,\cdot)\|^2_{L^2(Q;\R^M)}
  +\|f^\eta_n(x,\cdot)\|_{L^2(Q;\R^M)}\|\nabla_y f^\eta_n(x,\cdot)\|_{L^2(Q;\R^M)} \right] \\
&\hspace{2cm}\leq \sqrt{2}\left[2\eta^2+\eta\|\nabla_y \unf u_n(x,\cdot)\|_{L^2(Q;\R^M)} \right],
\end{align*}
where in the last step we used \Cref{lem:truncation}, \eqref{eq:est_v_eta_n_k}, and the definition of $k_n$. 
This concludes the proof of the slicing theorem.
\end{proof}


We now continue with the proof of \Cref{thm:liminf}.

\textbf{Step 3: Passing to Limit.}
Using \Cref{thm:sequence_kxn}, we get that, for all $\eta>0$, it holds
\[
\liminf_{n\to\infty} \mathcal{E}_n(u_n) \geq \liminfn\int_\o{2\left[\int_Q W(y,u_n+v_n^\eta)\;dy\right]^{\frac{1}{2}}|\nabla u_n|\;dx}.
\]
We now want to pass to the limit in the above inequality as $n \to \infty$, replacing  $W$ by $W^\eta$. Observe that $v_n^\eta(x,\cdot)$ might not satisfy, for a.e $x\in\Omega$, the required bound to be an admissible competitor for the minimization problem defining $W^\eta$. Thus, we reason as follows.
Define
\[
\o_n^\eta \coloneqq \{x \in \o:\; \|\nabla_y \unf u_n(x,\cdot)\|_{L^2(Q;\R^M)} \leq \eta\}.
\]
Then, for each $x\in \o_n^\eta$, we have that $v_n^\eta(x,\cdot)\in\mathcal{A}_\eta$ (see Definition \ref{def:W_eta}).
We can estimate
\begin{align*}
\liminfn\int_\o{2\left[\int_Q W(y,u_n+v_n^\eta)\;dy\right]^{\frac{1}{2}}|\nabla u_n|\;dx}
	\geq \liminfn\int_{\o_n^\eta}{2\sqrt{W^\eta(u_n)}|\nabla u_n|\;dx} \\
 \geq \liminfn\int_{\o}{2\sqrt{W^\eta(u_n)}|\nabla u_n|\;dx} - \limsupn \int_{\o\setminus\o_n^\eta}{2\sqrt{W^\eta(u_n)}|\nabla u_n|\;dx}.
\end{align*}
We claim that
\[
\limsupn \int_{\o\setminus\o_n^\eta}{2\sqrt{W^\eta(u_n)}|\nabla u_n|\;dx}=0.
\]
Indeed, by Chebyshev's inequality we have that $$\mathcal{L}^N(\o\setminus\o_n^\eta) \leq \frac{1}{\eta^2}\int_\o\|\nabla_y \unf u_n(x,\cdot)\|^2_{L^2(Q;\R^M)}\;dx, $$
and thus, by \eqref{eq:grad_est_two_scale} we get
\begin{equation}\label{eq:meas_oeta_n}
\mathcal{L}^N(\o\setminus\o_n^\eta) \leq \frac{\delta^2_n}{\eta^2\e_n}.
\end{equation}
Using the upper bound on $u_n$ (see \eqref{eq:assumptions_un}), we obtain
$$\int_{\o\setminus\o_n^\eta}{2\sqrt{W^\eta(u_n)}|\nabla u_n|\;dx} \leq C_{M,\eta}\int_{\o\setminus\o_n^\eta}{|\nabla u_n|\;dx},$$
and applying H\"{o}lder and \eqref{eq:meas_oeta_n} we infer from \eqref{eq:grad_est} that
\[
\int_{\o\setminus\o_n^\eta}{2\sqrt{W^\eta(u_n)}|\nabla u_n|\;dx}\leq C_{M,\eta}\mathcal{L}^N(\o\setminus\o_n^\eta)^{\frac{1}{2}}\|\nabla u_n\|_{L^2(Q;\R^M)} \leq C_{M,\eta}\frac{\delta_n}{\e_n},
\]
which proves the claim because $\delta_n \ll \e_n$.

Therefore, for all $\eta>0$ we have that
\begin{equation}\label{eq:final_liminf}
   \liminfn \mathcal{E}_n (u_n) \geq \liminfn\int_{\o}{2\sqrt{W^\eta(u_n)}|\nabla u_n|\;dx}. 
\end{equation}

We specifically note that up to this point, we have not used the convergence of $\{u_n\}_n$ in any way. Now we do so, and by virtue of "classical" results on the Modica-Mortola functional (see \cite[Theorem 3.4]{FonsecaTartar}), we get that
\begin{equation}\label{eq:liminf_eta}
\limn \mathcal{E}_n(u_n) \geq  \sigma_\eta \text{Per}(\{u = a\}),
\end{equation}
where
\[
\sigma_\eta:= \min\left\{\int_{-1}^1 2\sqrt{{W}^{\eta}(\gamma)}|\gamma'|dt \,:\,
    \gamma\in \mathrm{Lip}_{\mathcal{Z}}([-1,1];\R^M), \gamma(-1)=a,\, \gamma(1)=b  \right\}.
\]


\textbf{Step 3: Concluding the Arguments}
We send $\eta \to 0$ in \eqref{eq:liminf_eta}.
Recall that
\[
\sigma_{\mathrm{hom}}:= \min\left\{\int_{-1}^1 2\sqrt{{W}_{\mathrm{hom}}(\gamma)}|\gamma'|dt \,:\,
    \gamma\in \mathrm{Lip}_{\mathcal{Z}}([-1,1];\R^M), \gamma(-1)=a,\, \gamma(1)=b  \right\}.
\]
We claim the following.

\begin{proposition}
Let
\[
\sigma_0\coloneqq \lim_{\eta \to 0} \sigma_\eta = \sup_{\eta>0}{\sigma_\eta}.
\]
Then $\sigma_0 = \sigma_{\mathrm{hom}}$.
\end{proposition}

\begin{proof}
Note that $\sigma_\eta$ is increasing, and thus the limit exists.

By \Cref{thm:prop_W_eta}, for every $\eta>0$, $\sigma_\eta\leq \sigma_{\mathrm{hom}}$, and
this establishes the inequality $\sigma_0\leq \sigma_{\mathrm{hom}}$.

We prove the converse inequality. \Cref{prop:d0_character} gives us that $\sigma_0$ is the distance between $a$ and $b$ in a certain metric $\mathrm{d}_0$, and that there exists $\gamma_0\in  C^0([-1,1];\R^M)$  such that
\begin{equation}\label{eq:sigma0}
\sigma_0 = \sup_{\eta > 0} \int_{-1}^1 {2\sqrt{W^\eta(\gamma_0)}|\gamma_0'|dt}.
\end{equation}
Without loss of generality, due to parameterization invariance we can assume that the $\gamma_0(t)=a$ if and only if $t=-1$, and that $\gamma_0(t)=b$ if and only if $t=1$.
For $j \in \N\setminus\{0\}$, define
\[
T^a_j\coloneqq \{\, t\in[-1,1] \,:\, \gamma_0(t)\not\in \overline{B(a,1/j)}  \,\},
\]
\[
T^b_j\coloneqq \{\, t\in[-1,1] \,:\, \gamma_0(t)\not\in \overline{B(b,1/j)}  \,\},
\]
and
\[
T_j\coloneqq [-1,1]\setminus (T^a_j \cup T^b_j).
\]

Let
\[
L_j\coloneqq \int_{T^j} |\gamma'_0(t)|\, dt < \infty.
\]
By \Cref{thm:prop_W_eta}, we have uniform convergence of $W^\eta$ ot $W_{\mathrm{hom}}$ on compact sets.
Let $(\eta_j)_j$ be a decreasing sequence such that
\[
\|\sqrt{{W}^{\eta_j}}-\sqrt{W_{\mathrm{hom}}}\|_{C^0(K)} \leq \frac{1}{2j L_j}.
\]
for all $\eta<\eta_j$, where $K\subset\R^M$ is a compact set such that $\gamma_0(t)\in K$ for all $t\in[-1,1]$.
We observe that
\begin{align*}
     \sigma_0 &\geq \int_{-1}^1 {2\sqrt{W^{\eta_j}(\gamma_0)}|\gamma_0'|dt} \\
     & \geq \int_{T^j} {2\sqrt{W^{\eta_j}(\gamma_0)}|\gamma_0'|dt} \\
     &\geq \int_{T^j} {2\sqrt{W_{\mathrm{hom}}(\gamma_0)}|\gamma_0'|dt}
        -2\int_{T^j} \left|\sqrt{W^{\eta_j}(\gamma_0)}-\sqrt{{W}_{\mathrm{hom}}(\gamma_0)}\right||\gamma_0'|dt\\
     &\geq \int_{T^j}{2\sqrt{W_{\mathrm{hom}}(\gamma_0)}|\gamma_0'|dt}-\frac{1}{j}
\end{align*}
Taking $j \to \infty$ and using Monotone Convergence Theorem together with \eqref{eq:sigma0}, we conclude. 
\end{proof}


\section{Compactness}

In this section prove \Cref{thm:main_compactness}.

\begin{proof}[Proof of Theorem \ref{thm:main_compactness}]
Let $\{u_n\}_n \subset W^{1,2}(\Omega;\R^M)$ be such that
\[
\sup_{n\in\N} \mathcal{E}_n(u_n) = C<\infty.
\]

By \eqref{eq:final_liminf}, we obtain the uniform bound without using any information about the convergence:
$$\sup_n \int_{\o}{2\sqrt{W^\eta(u_n)}|\nabla u_n|\;dx}\leq C$$ 
Since $W^\eta$ has two wells $a,b$, we can apply arguments in \cite{FonsecaTartar} to extract a subsequence such that $u_n \to u \in BV(\o;\{a,b\})$ strongly in $L^1$.

\end{proof}


\section{Limsup inequality}

The main result of this section is the construction of a recovery sequence.

\begin{theorem}\label{prop:limsup}
Given $u\in BV(\Omega;\{a,b\})$, there exists a sequence $\{u_n\}_n\subset W^{1,2}(\Omega;\R^M)$ with $u_n\to u$ strongly in $L^1(\o;\R^M)$, such that
\[
\limsup_{n\to\infty} \mathcal{E}_n(u_n) \leq \mathcal{E}_\infty(u).
\]
\end{theorem}

The proof of the limsup inequality requires two technical results that are well known to the experts. For the reader's convenience, we state them in here. The first is an approximation result, stating that $C^2$ sets are dense both in configuration and in energy. This result is contained in the proof of \cite[Lemma 3.1]{baldo}.

\begin{proposition}\label{prop:approx_C2_sets}
Let $E\subset \Omega$ be a set with finite perimeter. Then, there exists a sequence of sets $\{E_n\}_n$, where each $E_n\subset\Omega$ satisfies
\begin{itemize}
\item $\partial E_n\cap \Omega$ is of class $C^2$;
\item $\hno(\partial\Omega\cap \partial E_n)=0$,
\end{itemize}
such that $E_n\to E$ with respect to the $L^1$ topology, and
\[
\lim_{n\to\infty} \mathcal{E}_\infty(u_n) = \mathcal{E}_\infty(u),
\]
where $u_n\coloneqq \mathbbmss{1}_{E_n}$, and $u\coloneqq \mathbbmss{1}_E$.
\end{proposition}

The next result ensures that, up to a small error, it is possible to reparametrize a curve in such a way that the energy functional is bounded by the limiting energy. This was originally used in the article by Modica (see \cite[Proof of Proposition 2]{modica87}). Here we give the version used in \cite[Lemma 4.5]{CriGra}, since it clearly states the estimates that we will need. Note that none of the assumptions on $W$ required in \cite{CriGra} are actually used, other than continuity in the second variable. Moreover, the lower bound for $\tau$ follows easily from the definition of $\tau$ given in the proof of the result.

\begin{lemma}\label{lem:reparametrization}
Fix $\lambda>0$, $\varepsilon>0$.
Let $\gamma\in C^1([-1,1];\R^M)$, with $\gamma(-1)=a$, $\gamma(1)=b$, and  $\gamma'(s) \neq 0$ for all $s \in (-1,1)$.
Then, there exist $\tau>0$, and $C>0$, with
\[
C\varepsilon \leq \tau \le \frac{\e}{\sqrt{\lambda}}\int_{-1}^1 |\gamma'(t)|\,dt,
\]
and $g \in C^1((-\tau,\tau); [-1,1])$ such that
\[
(g'(t))^2 = \frac{\lambda + W(x,\gamma(g(t)))}{\e^2 |\gamma'(g(t))|^2}
\]
for all $t\in(-\tau,\tau)$, $g(-\tau)=-1$, $g(\tau)=1$, and
\begin{multline*}
\int_{-\tau}^{\tau} \left[\frac{1}{\e} W_{\mathrm{hom}}\left(\gamma(g(t)) \right) +\e |\gamma'\left( g(t) \right)|^2\left(g'(t)\right)^2 \right] \,dt \\ 
\le \int_{-1}^1 2 \sqrt{W_{\mathrm{hom}}(\gamma(s))}|\gamma'(s)|\,ds + 2\sqrt{\lambda} \int_{-1}^1 |\gamma'(s)| \,ds.
\end{multline*}
\end{lemma}

We are now in position to prove the existence of a recovery sequence.

\begin{proof}[Proof of \Cref{prop:limsup}]
Let $A\coloneqq \{u=a\}$.

\begin{figure}
\includegraphics[scale=0.7]{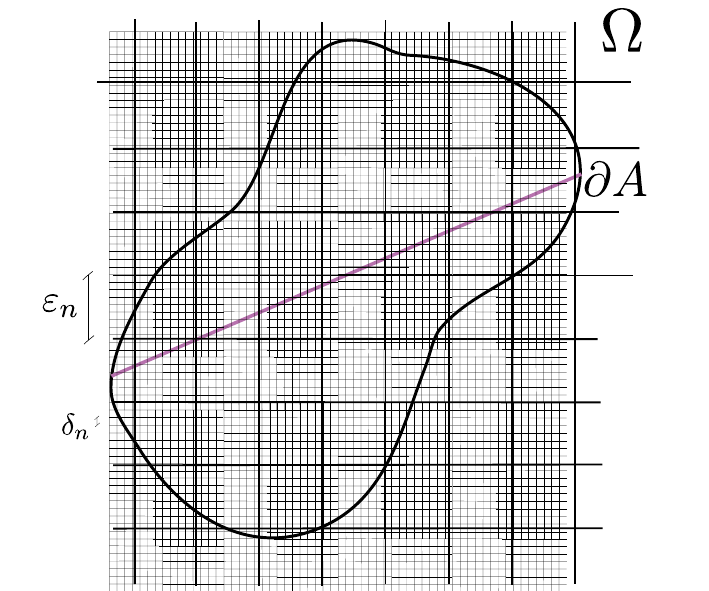}
\caption{The two different scales $\varepsilon_n$ and $\delta_n$.}
\label{fig:microstructure}
\end{figure}

\textbf{Step 1.} Using Proposition \ref{prop:approx_C2_sets} and a diagonalization argument, we note that it suffices to prove the result for $u\in BV(\Omega;\{a,b\})$ such that  $\partial A \cap \Omega$ is of class $C^2$, and $\hno(\partial A \cap \partial \o) = 0$.\\

\textbf{Step 2.} Let $u\in BV(\Omega;\{a,b\})$ be as in Step 1. 
Fix $\eta>0$.
Let $\gamma\in C^1([-1,1];\R^N)$ with $\gamma(-1)=a$, and $\gamma(1)=b$, be such that
\begin{equation}\label{eq:gamma_sigma_hom}
\int_{-1}^1 2\sqrt{{W}_{\mathrm{hom}}(\gamma(t))}|\gamma'(t)|dt \leq \sigma_{\mathrm{hom}} + \eta.
\end{equation}
Without loss of generality, we can assume that $\gamma'(s)\neq0$ for all $s\in(-1,1)$.
For each $n\in\N$, let $\tau_n>0$, and $g_n\in C^1([0,\tau_n];\R^M)$ be those given by Lemma \ref{lem:reparametrization} relative to the choice of
\begin{equation}\label{eq:eps_lambda}
\varepsilon := \varepsilon_n,\quad\quad\quad
\lambda:= \left( \frac{\eta}{L(\gamma)} \right)^2,
\end{equation}
where
\[
L(\gamma)\coloneqq \int_{-1}^1 |\gamma'(s)|\, ds < \infty.
\]
Let $\mathrm{dist}(\cdot, \partial A):\R^N\to\R$ be the signed distance function from $\partial A$. Note that $\mathrm{dist}(\cdot,\partial A)$ is of class $C^1$, since $\partial A$ is of class $C^2$.
For $n\in\N$, define
\begin{equation}\label{eqn:rec_seq}
u_n(x) \coloneqq \begin{dcases}
	b & \mathrm{dist}(x,\partial A)>\tau_n,\\
	\gamma\left(g_n(\mathrm{dist}(x,\partial A))\right) &
        |\mathrm{dist}(x,\partial A)|\leq \tau_n, \\
	a & \mathrm{dist}(x,\partial A)<-\tau_n.
\end{dcases}
\end{equation}
We claim that the sequence $\{u_n\}_n$ satisfies the required properties, up to an error $\eta$. Observe that each $u_n\in W^{1,2}(\Omega;\R^M)$. Moreover, by using the fact that $\tau_n\to0$ as $n\to\infty$, we have that $u_n\to u$ strongly in $L^1(\o;\R^M)$.

We prove the convergence of the energies.
Define
\[
A_n \coloneqq \{x \in \Omega \,:\, |\mathrm{dist}(x,\partial A)|\leq \tau_n\},
\]
and note that
\begin{equation}\label{eq:An}
|A_n|\leq \hno(\partial A)\tau_n.
\end{equation}
For each $n\in\N$, consider the set of indexes
\[
I_n \coloneqq \{ i\in\N \,:\, z_i+\delta_n Q \subset A_n,\, z_i \in \delta_n \Z^N \},
\]
\[
B_n \coloneqq \{ i\in\N \,:\, z_i+\delta_n Q \cap A_n \neq \emptyset,\, z_i \in \delta_n \Z^N \}\setminus I_n.
\]
Let $K_n\coloneqq \# I_n$, and observe that
\begin{equation}\label{eq:Kn}
\lim_{n\to\infty} K_n\left(\left\lfloor\frac{|A_n|}{\delta_n^N}\right\rfloor\right)^{-1} = 1.
\end{equation}
Write
\[
A_n : = \bigcup_{i\in I_n} \left(z_i+\delta_n Q\right) \cup R_n,
\]
where
\[
R_n\coloneqq \bigcup_{j\in B_n} \left(z_j+\delta_n Q\right)\cap A_n.
\]
Then,
\begin{equation}\label{eqn:Remainder}
	|R_n|\leq C_{A}\delta_n,
\end{equation} 
where $C_A>0$ is a constant depending on $\partial A$.
\vspace{0.3cm}

\textbf{Step 3.} We claim that there exists a dimensional constant $C_N>0$ such that
\begin{equation}\label{eq:est_limsup_step3}
|u_n(x)-u_n(z_i)|\leq \omega_\gamma\left(C_N\frac{\delta_n}{\e_n}\right),
\end{equation}
for all $x \in z_i+\delta_n Q$, and all $n\in\N$, and $i\in I_n$, where $\omega_\gamma:[0,\infty)\to[0,\infty)$ is the modulus of continuity of $\gamma$.
Indeed,
\begin{align*}
|u_n(x)-u_n(z_i)|&= \left|\gamma\left(g_n(\mathrm{dist}(x,\partial A))\right) - 
    \gamma\left(g_n(\mathrm{dist}(z_i,\partial A))\right)\right| \\
&\leq \omega_\gamma\left(\frac{1}{\e_n}(|\mathrm{dist}(x,\partial A) - \mathrm{dist}(z_i,\partial A)|)\right) \\
&\leq \omega_\gamma\left(\frac{1}{\e_n}(|x-z_i|)\right) \\
&\leq \omega_\gamma\left(C_N\frac{\delta_n}{\e_n}\right),
\end{align*}
where in the second step we used the fact that $|g'_n|\leq C/\varepsilon_n$.
\vspace{0.3cm}

\textbf{Step 4.} We claim that
\[
\lim_{n\to\infty}\left|\frac{1}{\e_n}\int_{A_n} \left[ W\left(\frac{x}{\delta_n},u_n(x)\right)
    - W_{\mathrm{hom}}(u_n(x)) \right] \;dx\right|= 0.
\]
By using the unfolding operator restricted to $A_n$, we get
\begin{align*}
\int_{A_n}{W\left(\frac{x}{\delta_n},u_n(x)\right)\;dx}
&= \sum_{i=1}^{K_n}\int_{z_i+\delta_n Q}\int_Q W(y,\unf u_n)\;dydx
    + \int_{R_n}{W\left(\frac{x}{\delta_n},u_n(x)\right)\;dx} \nonumber \\
&= \delta_n^N\sum_{i=1}^{K_n} \int_Q W(y,u_n(z_i+\delta_n y))\;dy\
    + \int_{R_n}{W\left(\frac{x}{\delta_n},u_n(x)\right)\;dx}.
\end{align*}
Thus, we can write
\begin{align}\label{eqn:splitting}
&\int_{A_n}{ \left[ W\left(\frac{x}{\delta_n},u_n(x)\right) - W_{\mathrm{hom}}(u_n(x)) \right]  \;dx} \nonumber \\
&\hspace{2cm}= \delta_n^N\sum_{i=1}^{K_n} \left[\int_Q W(y,u_n(z_i+\delta_n y))\;dy - W_{\mathrm{hom}}(u_n(z_i))\right]\nonumber\\
&\hspace{4cm}+ \sum_{i=1}^{K_n}\left[\delta_n^N W_{\mathrm{hom}}(u_n(z_i))
    - \int_{z_i+\delta_n Q}{W_{\mathrm{hom}}(u_n(x))\;dx}\right]\nonumber\\
&\hspace{4cm}+\int_{R_n}{\left[W\left(\frac{x}{\delta_n},u_n(x)\right)-W_{\mathrm{hom}}(u_n(x))\right]\;dx}\nonumber\\
&\hspace{2cm}=: J_n^1 + J_n^2 + J_n^3.
\end{align}
We now estimate the three terms $J_n^1, J_n^2$, and $J_n^3$ separately.
Since $\gamma$ is continuous on a compact set, we note that by definition of the recovery sequence
\begin{equation}\label{eq:un_unif_bounded}
\sup_{n\in\N}\|u_n\|_{L^\infty} \leq M < \infty.
\end{equation}
This implies that
\begin{equation}\label{eqn:J3}
	|J_n^3|\leq C|R_n|\leq C\delta_n,
\end{equation}
where in the last inequality we used \eqref{eqn:Remainder}.\\

Now we estimate $J_n^2$. Let $\omega_{\mathrm{hom}}$ be a modulus of continuity of $W_{\mathrm{hom}}$ in $B(0,M)$, where $M>0$ is the constant in \eqref{eq:un_unif_bounded}. Then,
\begin{align}\label{eqn:J2}
|J_n^2| &= \left|\int_{\bigcup_{i=1}^{K_n} z_i+\delta_n Q}{W_{\mathrm{hom}}(u_n(z_i))-W_{\mathrm{hom}}(u_n(x))\;dx}\right| \nonumber \\
&\leq \int_{\bigcup_{i=1}^{K_n} z_i+\delta_n Q}{\omega_{\mathrm{hom}}(|u_n(z_i)-u_n(x)|)\;dx} \nonumber \\
& \leq K_n\delta_n^N\omega_{\mathrm{hom}}\left(\omega_\gamma\left(C_N\frac{\delta_n}{\e_n}\right)\right) \nonumber \\
 & \leq  C\e_n\omega_{\mathrm{hom}}\left(\omega_\gamma\left(C_N\frac{\delta_n}{\e_n}\right)\right),
\end{align} 
where in the previous to last inequality we used Step 3, while in the last inequality we used \eqref{eq:An} together with \eqref{eq:Kn}.

Finally, we estimate $J_n^1$. Using the definition of $W_{\mathrm{hom}}$, we write
\[
\int_Q W(y,u_n(z_i+\delta_n y))\;dy- W_{\mathrm{hom}}(u_n(z_i))
    =\int_Q \left[ W(y,u_n(z_i+\delta_n y))-W(y,u_n(z_i)) \right]\;dy.
\]
Thus, using a similar argument to that in \eqref{eqn:J2}, we obtain
\begin{equation}\label{eqn:J1}
	|J_n^1|\leq C\e_n\omega_{W}\left(\omega_\gamma\left(C_N\frac{\delta_n}{\e_n}\right)\right),
\end{equation}
where $\omega_W$ is a modulus of continuity of $W$ in $B(0,M)$, where $M>0$ is the constant in \eqref{eq:un_unif_bounded}.

Combining \eqref{eqn:J1}\eqref{eqn:J2}\eqref{eqn:J3}, we get
\begin{align}\label{eq:estimate_I123}
&\left|\frac{1}{\e_n}\int_{A_n}{W\left(\frac{x}{\delta_n},u_n(x)\right)- W_{\mathrm{hom}}(u_n(x))  \;dx}\right|
    \nonumber \\
&\hspace{2cm}\leq \frac{1}{\e_n}\left(|J_n^1|+|J_n^2|+|J_n^3|\right) \nonumber \\
&\hspace{2cm}\leq C\left(\omega_{W}\left(\omega_\gamma\left(C\frac{\delta_n}{\e_n}\right)\right)
    +\omega_{\mathrm{hom}}\left(\omega_\gamma\left(C\frac{\delta_n}{\e_n}\right)\right)
    + \frac{\delta_n}{\e_n} \right) \nonumber \\
&\hspace{2cm}\to 0,
\end{align}
as $n\to\infty$, where in the last step we used the fact that $\delta_n\ll\varepsilon_n$, together with $\lim_{t\to 0}\omega_W(t) = \lim_{t\to 0}\omega_{\mathrm{hom}}(t) = 0$.\\

\textbf{Step 5.} Using Step 4, and the coarea formula (see \cite[Theorem 2.93 and Remark 2.94]{AFP}), we get that
\begin{align*}
\limsupn \mathcal{E}_n(u_n) &=
    \limsupn \int_\o \left[ \frac{1}{\e_n}W\left(\frac{x}{\delta_n},u_n\right)
        +\e_n|\nabla u_n|^2 \right]\;dx \\
&\leq \limsupn\int_\o\left[\frac{1}{\e_n}W_{\mathrm{hom}}(u_n)+\e_n|\nabla u_n|^2 \right]\;dx \\
&\hspace{2cm}+ \limsupn\left| \frac{1}{\e_n}\int_{\o}
    \left[ W\left(\frac{x}{\delta_n},u_n(x)\right)
    - W_{\mathrm{hom}}(u_n(x)) \right] \;dx\right| \\
&= \limsupn\int_\o\left[\frac{1}{\e_n}W_{\mathrm{hom}}(u_n)+\e_n|\nabla u_n|^2 \right]\;dx \\
&=\limsupn \int_{-\tau_n}^{\tau_n} \left[\frac{1}{\e_n}W_{\mathrm{hom}}(\gamma(g_n(s))) + \e_n|\gamma'(g_n(s)) g'_n(s)|^2 \right]\cdot \\
&\hspace{2cm}\cdot \mathcal{H}^{N-1}(\{ x\in \Omega : \mathrm{dist}(x,\partial A) = s \}) \;ds \\
&\leq \limsupn \sup_{|s|\leq \tau_n} \mathcal{H}^{N-1}(\{ x\in \Omega : \mathrm{dist}(x,\partial A) = s \})  \cdot \\
&\hspace{2cm} \cdot\int_{-\tau_n}^{\tau_n} \left[\frac{1}{\e_n}W_{\mathrm{hom}}(\gamma(g_n(s))) + \e_n|\gamma'(g_n(s)) g'_n(s)|^2 \right]\, ds\\
&\leq \limsupn \sup_{|s|\leq \tau_n} \mathcal{H}^{N-1}(\{ x\in \Omega : \mathrm{dist}(x,\partial A) = s \})  \cdot \\
&\hspace{2cm} \cdot \left[ \int_{-1}^1 2\sqrt{{W}_{\mathrm{hom}}(\gamma(t))}|\gamma'(t)|dt + 2\sqrt{\lambda} L(\gamma)  \right] \\
&\leq \mathcal{H}^{N-1}(\partial A)[\sigma_{\mathrm{hom}} + 3\eta],
\end{align*}
where the last inequality follows from \eqref{eq:estimate_I123}, together with \eqref{eq:eps_lambda}, Lemma \ref{lem:reparametrization}, and the fact that, since $\partial A$ is of class $C^2$, it holds
\[
\lim_{s\to 0} \mathcal{H}^{N-1}(\{ x\in \Omega : \mathrm{dist}(x,\partial A) = s \}) = \mathcal{H}^{N-1}(\partial A).
\]
Since $\eta>0$ is arbitrary, we conclude.
\end{proof}


\section{The mass constrained functional}

As it is usually the case, the strategy of the proof for the Gamma-limit in the unconstrained case is stable enough to be used for the mass constrained functional, with minor changes. The liminf inequality follows from exactly the same proof. What needs to be checked is that, in the construction of the recovery sequence, it is possible to adjust the mass of the sets constructed in the proof of Proposition \ref{prop:limsup}.
There are some standard ways to do that, which can be found in the paper by Fonseca and Tartar (see \cite{FonsecaTartar}), and in the paper by Ishige (see \cite{Ishige}, see also \cite{baldo}). However, such strategies are heavily based on the assumption that potential $W$ enjoys regularity that we do not require in this paper. Here, we use a different argument that does not require such assumption. 

\begin{lemma}
Fix $m\in(0,|\Omega|)$. Let $u\in BV(\Omega;\{a,b\})$ with $|\{u=a\}|=m$. Then, there exists a sequence $\{u_n\}_n\subset W^{1,2}(\Omega;\R^M)$ with $u_n\to u$ strongly in $L^1(\o;\R^M)$ such that
\[
\limsup_{n\to\infty} 
\widehat{\mathcal{E}}_n(u_n) \leq \widehat{\mathcal{E}}_\infty(u),
\]
where $\widehat{\mathcal{E}}_n$ and $\widehat{\mathcal{E}}_\infty$ are defined in Definition \ref{def:mass_constrain_functional}.
\end{lemma}

\begin{proof}
We show how to modify Proposition \ref{prop:approx_C2_sets} and the definition of the function $u_n$ in Step 2 of the proof of Proposition \ref{prop:limsup} in order to get the recovery sequence satisfying the mass constraint.\\

\textbf{Step 1.} Let $E\coloneqq\{u=a\}$. Let $\{E_n\}_n$ be the sequence provided by Proposition \ref{prop:approx_C2_sets}. We would like to modify this in such a way that the required mass constraint is satisfied. The strategy we use is a variant of an idea by Ryan Murray.

Let $x_0, x_1\in \Omega$ be points of density one and zero for $E$, respectively. i.e.,
\[
\lim_{r\to0} \frac{|E\cap B(x_1,r)|}{|B(x_1,r)|} = 1\quad\quad\text{and}\quad\quad
\lim_{r\to0} \frac{|E\cap B(x_0,r)|}{|B(x_0,r)|} = 0.
\]
Then,
there exists $R>0$ such that
\begin{equation}\label{eq:density_estimates}
\frac{3}{4}\leq  \frac{|E\cap B(x_1,r)|}{|B(x_1,r)|} \leq 1,\quad\quad\quad\quad
0\leq \frac{|E\cap B(x_0,r)|}{|B(x_0,r)|} \leq \frac{1}{4},
\end{equation}
for all $r\leq R$. Without loss of generality, up to decreasing the value of $R>0$, we can assume that $B(x_1,r)\cap B(x_0,r)=\emptyset$. For $r\leq R$, let
\[
\widetilde{E}_r\coloneqq E\cup B(x_1,r) \setminus B(x_0,r).
\]
Note that
\[
\lim_{r\to0}|P(E_r; \Omega) - P(E;\Omega)|=0,\quad\quad\quad
\| \mathbbmss{1}_{E_r} - \mathbbmss{1}_{E} \|_{L^1(\Omega)} = 0.
\]
Let $\{F^r_n\}_n$ be the sequence of sets given by Proposition \ref{prop:approx_C2_sets} relative to $E_r$. Then, for $n$ large, we get that
\[
|P(F^r_n; \Omega) - P(E_r;\Omega)|<r,\quad\quad\quad
\| \mathbbmss{1}_{E^r_n} - \mathbbmss{1}_{E^r} \|_{L^1(\Omega)} < r.
\]
Since such sets $F^r_n$ are obtained with the standard procedure of mollification of the characteristic function of $E_r$, and by taking a super level set, we get that there exists $\bar{n}\in\N$ such that
\[
B\left(x_0, \left(\frac{4}{5}\right)^N r \right) \subset \Omega\setminus F^r_n\quad\quad\text{and}\quad\quad
B\left(x_1, \left(\frac{4}{5}\right)^N r \right) \subset F^r_n,
\]
for all $n\geq\bar{n}$. Now, assume that $|F^r_n|<m$. Let $s_n>0$ be such that
$|F^r_n| + |B(0,s_n)| = m$. We claim that
\[
s_n < \left(\frac{4}{5}\right)^N r,
\]
for all $n\geq\bar{n}$. Indeed, for $n$ large enough, using \eqref{eq:density_estimates} it holds that
\[
|E_r| - |E| < \frac{3}{4}|B(x_1,r)|.
\]
Therefore, considering the set
\[
\widetilde{F}^r_n\coloneqq F^r_n \cup B(x_0, s_n),
\]
we get that $|\widetilde{F}^r_n|=m$, for all $n\geq \bar{n}$.
Using a diagonal argument, we obtain the desired conclusion.
A similar argument is used to fix the mass in the case where $|F^r_n|>m$. This sequence satisfies the required properties.\\

\textbf{Step 2.} It can be shown that, for each $n\in\N$, there exists $v_n\in\R$ such that the function
\begin{equation*}
\widetilde{u}_n(x) \coloneqq \begin{dcases}
	b & \mathrm{dist}(x,\partial A)>\tau_n,\\
	\gamma\left(g_n(\mathrm{dist}(x,\partial A) + v_n)\right) &
        |\mathrm{dist}(x,\partial A)|\leq \tau_n, \\
	a & \mathrm{dist}(x,\partial A)<-\tau_n.
\end{dcases}
\end{equation*}
is such that 
\[
\int_\Omega u_n(x)\, dx = ma + (1-m)b.
\]
Using the fact that $v_n\to 0$ as $n\to\infty$, it is possible to check that all of the steps in the proof of Proposition \ref{prop:limsup} can be carried out in a similar way. This allows to conclude.
\end{proof}


\section*{Acknowledgements}
I.F. was partially supported under NSF-DMS1906238 and NSF-DMS2205627. \\L.G. was funded by the Deutsche Forschungsgemeinschaft–320021702/GRK2326– Energy, Entropy, and Dissipative Dynamics (EDDy). \\

\textbf{Data Availability}: Data sharing not applicable to this article as no datasets were generated or analysed during the current study.\\

\textbf{Ethics Statement}: The authors have no conflicts of interest to declare relevant to this article.\\


\bibliographystyle{siam}
\bibliography{Supercritical}

\begin{thebibliography}{10}

\bibitem{AFP}
{\sc L.~Ambrosio, N.~Fusco, and D.~Pallara}, {\em Functions of bounded
  variation and free discontinuity problems}, Oxford Mathematical Monographs,
  The Clarendon Press, Oxford University Press, New York, 2000.

\bibitem{AnsBraChi1}
{\sc N.~Ansini, A.~Braides, and V.~Chiad\`o~Piat}, {\em Interactions between
  homogenization and phase-transition processes}, Tr. Mat. Inst. Steklova, 236
  (2002), pp.~386--398.

\bibitem{AnsBraChi2}
{\sc N.~Ansini, A.~Braides, and V.~Chiad\`o~Piat}, {\em Gradient theory of
  phase transitions in composite media}, Proc. Roy. Soc. Edinburgh Sect. A, 133
  (2003), pp.~265--296.

\bibitem{anzellotti1993asymptotic}
{\sc G.~Anzellotti and S.~Baldo}, {\em Asymptotic development by
  {$\Gamma$}-convergence}, Appl. Math. Optim., 27 (1993), pp.~105--123.

\bibitem{BacMarZep22}
{\sc A.~Bach, T.~Esposito, R.~Marziani, and C.~I. Zeppieri}, {\em Gradient
  damage models for heterogeneous materials}, 2022.
\newblock cvgmt preprint.

\bibitem{Bach2022}
\leavevmode\vrule height 2pt depth -1.6pt width 23pt, {\em Interaction Between
  Oscillations and Singular Perturbations in a One-Dimensional Phase-Field
  Model}, Springer International Publishing, Cham, 2022, pp.~3--31.

\bibitem{baldo}
{\sc S.~Baldo}, {\em Minimal interface criterion for phase transitions in
  mixtures of {C}ahn-{H}illiard fluids}, Ann. Inst. H. Poincar\'e Anal. Non
  Lin\'eaire, 7 (1990), pp.~67--90.

\bibitem{Bou}
{\sc G.~Bouchitt\'e}, {\em Singular perturbations of variational problems
  arising from a two-phase transition model}, Appl. Math. Optim., 21 (1990),
  pp.~289--314.

\bibitem{Braides}
{\sc A.~Braides}, {\em {$\Gamma$}-convergence for beginners}, vol.~22 of Oxford
  Lecture Series in Mathematics and its Applications, Oxford University Press,
  Oxford, 2002.

\bibitem{braides2008asymptotic}
{\sc A.~Braides and L.~Truskinovsky}, {\em Asymptotic expansions by
  {$\Gamma$}-convergence}, Contin. Mech. Thermodyn., 20 (2008), pp.~21--62.

\bibitem{BradiesZeppieri}
{\sc A.~Braides and C.~I. Zeppieri}, {\em Multiscale analysis of a prototypical
  model for the interaction between microstructure and surface energy},
  Interfaces Free Bound., 11 (2009), pp.~61--118.

\bibitem{CioDamGri08}
{\sc D.~Cioranescu, A.~Damlamian, and G.~Griso}, {\em The periodic unfolding
  method in homogenization}, SIAM J. Math. Anal., 40 (2008), pp.~1585--1620.

\bibitem{clyne2019introduction}
{\sc T.~W. Clyne and D.~Hull}, {\em An introduction to composite materials},
  Cambridge university press, 2019.

\bibitem{CristoferiFonsecaGanedi}
{\sc R.~Cristoferi, I.~Fonseca, and L.~Ganedi}, {\em Homogenization and phase
  separation with space dependent wells -- the subcritical case}, 2022.
\newblock https://arxiv.org/abs/2205.12893.

\bibitem{CriFonHagPop}
{\sc R.~Cristoferi, I.~Fonseca, A.~Hagerty, and C.~Popovici}, {\em A
  homogenization result in the gradient theory of phase transitions},
  Interfaces Free Bound., 21 (2019), pp.~367--408.

\bibitem{CriGra}
{\sc R.~Cristoferi and G.~Gravina}, {\em Sharp interface limit of a multi-phase
  transitions model under nonisothermal conditions}, Calc. Var. Partial
  Differential Equations, 60 (2021), pp.~Paper No. 142, 62.

\bibitem{Dalmasobook}
{\sc G.~Dal~Maso}, {\em An introduction to {$\Gamma$}-convergence}, vol.~8 of
  Progress in Nonlinear Differential Equations and their Applications,
  Birkh\"{a}user Boston, Inc., Boston, MA, 1993.

\bibitem{EG}
{\sc L.~C. Evans and R.~F. Gariepy}, {\em Measure theory and fine properties of
  functions}, Textbooks in Mathematics, CRC Press, Boca Raton, FL, revised~ed.,
  2015.

\bibitem{FonPopo}
{\sc I.~Fonseca and C.~Popovici}, {\em Coupled singular perturbations for phase
  transitions}, Asymptot. Anal., 44 (2005), pp.~299--325.

\bibitem{FonsecaTartar}
{\sc I.~Fonseca and L.~Tartar}, {\em The gradient theory of phase transitions
  for systems with two potential wells}, Proc. Roy. Soc. Edinburgh Sect. A, 111
  (1989), pp.~89--102.

\bibitem{Giusti}
{\sc E.~Giusti}, {\em Minimal surfaces and functions of bounded variation},
  vol.~80 of Monographs in Mathematics, Birkh\"{a}user Verlag, Basel, 1984.

\bibitem{Gurtin}
{\sc M.~E. Gurtin}, {\em Some results and conjectures in the gradient theory of
  phase transitions}, in Metastability and incompletely posed problems
  ({M}inneapolis, {M}inn., 1985), vol.~3 of IMA Vol. Math. Appl., Springer, New
  York, 1987, pp.~135--146.

\bibitem{Hagerty}
{\sc A.~Hagerty}, {\em A note on homogenization effects on phase transition
  problems}.
\newblock arXiv preprint arXiv:1811.07357, 2018.

\bibitem{Ishige}
{\sc K.~Ishige}, {\em Singular perturbations of variational problems of vector
  valued functions}, Nonlinear Anal., 23 (1994), pp.~1453--1466.

\bibitem{KohnStern}
{\sc R.~V. Kohn and P.~Sternberg}, {\em Local minimisers and singular
  perturbations}, Proc. Roy. Soc. Edinburgh Sect. A, 111 (1989), pp.~69--84.

\bibitem{MaggiBook}
{\sc F.~Maggi}, {\em Sets of finite perimeter and geometric variational
  problems}, vol.~135 of Cambridge Studies in Advanced Mathematics, Cambridge
  University Press, Cambridge, 2012.
\newblock An introduction to geometric measure theory.

\bibitem{Mar22}
{\sc R.~Marziani}, {\em {$\Gamma$}-convergence and stochastic homogenisation of
  phase-transition functionals}, 2022.
\newblock cvgmt preprint.

\bibitem{modica87}
{\sc L.~Modica}, {\em The gradient theory of phase transitions and the minimal
  interface criterion}, Arch. Rational Mech. Anal., 98 (1987), pp.~123--142.

\bibitem{modica77}
{\sc L.~Modica and S.~Mortola}, {\em Un esempio di {$\Gamma
  ^{-}$}-convergenza}, Boll. Un. Mat. Ital. B (5), 14 (1977), pp.~285--299.

\bibitem{Stern}
{\sc P.~Sternberg}, {\em The effect of a singular perturbation on nonconvex
  variational problems}, Arch. Ration. Mech. Anal., 101 (1988), pp.~209--260.

\bibitem{Vis06}
{\sc A.~Visintin}, {\em Towards a two-scale calculus}, ESAIM Control Optim.
  Calc. Var., 12 (2006), pp.~371--397.

\bibitem{Vis07}
{\sc A.~Visintin}, {\em Two-scale convergence of some integral functionals},
  Calc. Var. Partial Differential Equations, 29 (2007), pp.~239--265.

\bibitem{SternbergZuniga}
{\sc A.~Zuniga and P.~Sternberg}, {\em On the heteroclinic connection problem
  for multi-well gradient systems}, J. Differential Equations, 261 (2016),
  pp.~3987--4007.

\end{thebibliography}

\end{document}